\documentclass[12pt]{amsart}\usepackage{times}
\newcommand{\br}{\hfill\break}
\usepackage{amssymb}

\numberwithin{equation}{section}

\newcommand{\lab}[1]{\label{#1}}

\newtheorem{theorem}{Theorem}[section]
\newtheorem{lemma}[theorem]{Lemma}

\newenvironment{Proof}{\removelastskip\par\medskip
\noindent{\em Proof.} \rm}{\penalty-20\null\hfill$\square$\par\medbreak}

\newcommand{\di}{\mathrm{d}}

\newcommand{\R}{\mathbb{R}}

\newcommand{\g}{\mathfrak{g}}
\newcommand{\G}{\mathcal{G}}
\newcommand{\sbullet}{{\,\scriptstyle{{}^{{}_\bullet}}\,}}
\newcommand{\dslash}{{/\!/}}
\newcommand{\PM}{P\!M}
\newcommand{\etasm}{{\eta^{\mathrm{smooth}}}}
\newtheorem{Thm}[theorem]{Theorem}
\newtheorem{Prop}[theorem]{Proposition}
\newtheorem{Lem}[theorem]{Lemma}

\theoremstyle{remark}
\newtheorem{Rem}[theorem]{Remark}
\newtheorem*{Ack}{Acknowledgment}

\theoremstyle{definition}

\newtheorem{Exa}[theorem]{Example}

\newtheorem{Ass}[theorem]{Assumption}

\newcommand{\bbR}{{\mathbb{R}}}

\newcommand\Poiss[2]{\left\{{ {#1}\,,\,{#2} }\right\}}
\newcommand{\de}{\partial}
\newcommand{\pfrac}[1]{\genfrac(){}{}{#1}}

\newcommand{\bx}{{\mathbf{x}}}
\newcommand{\by}{{\mathbf{y}}}
\newcommand{\bv}{{\mathbf{v}}}
\newcommand{\ba}{{\mathbf{a}}}
\newcommand{\bb}{{\mathbf{b}}}
\newcommand{\bX}{{\mathbf{X}}}

\newcommand{\bfeta}{{\boldsymbol{\eta}}}
\newcommand{\bbeta}{{\boldsymbol{\beta}}}

\newcommand{\tT}{\widehat{T}}

\newcommand{\fg}{\mathfrak{g}}
\newcommand{\frP}{\mathfrak{P}}

\newcommand{\calG}{\mathcal{G}}
\newcommand{\tcalG}{\Tilde{\mathcal{G}}}
\newcommand{\ttcalG}{\Tilde{\tcalG}}
\newcommand{\ucalG}{{\mathcal{G}}}
\newcommand{\bucalG}{\boldsymbol{\calG}}

\newcommand{\matP}{\mathsf{P}}

\newcommand{\pepsilon}[3]{\epsilon^{{#1}{#2}}_{\phantom{{#1}{#2}}{#3}}}

\newcommand{\dd}{\mathrm{d}}
\newcommand{\T}{{T}}
\DeclareMathOperator{\ad}{ad}
\DeclareMathOperator{\Ad}{Ad}
\DeclareMathOperator{\hol}{Hol}

\begin{document}

\title[Poisson sigma models and symplectic groupoids]
{Poisson sigma models and symplectic groupoids}
\author[A. S. Cattaneo]{Alberto S. Cattaneo}
\thanks{A.~S.~C. acknowledges partial support of SNF Grant
No.~2100-055536.98/1}

\address{
A. S. C.: Mathematisches Institut\br
Universit\"at Z\"urich\br
CH-8057 Z\"urich, Switzerland}
\email{asc@math.unizh.ch}
\author[G. Felder]{Giovanni Felder}
\address{
G. F.:
D-MATH,\br
ETH-Zentrum\br
CH-8092 Z\"urich, Switzerland}

\email{felder@math.ethz.ch}

\begin{abstract}
We consider the Poisson sigma model associated to a Poisson
manifold. The perturbative
quantization of this model yields the Kontsevich
star product formula. We study here the classical model
in the Hamiltonian formalism. The phase space is
the space of leaves of a Hamiltonian foliation and 
has a natural groupoid structure. If it is a manifold then
it is a symplectic groupoid for the
given Poisson manifold. We study various families of examples.
In particular, a global symplectic groupoid for a
general class of two-dimensional Poisson domains is constructed.
\end{abstract}

\maketitle

\section{Introduction}\lab{S-1}
The notion of a symplectic groupoid \cite{Ka}, \cite{W} 
was introduced as part of a program to
quantize Poisson manifolds. It is modeled on the
following basic construction. 

Let $\g$ be a finite dimensional real
Lie algebra. Then its dual space $\g^*$ carries a Poisson
structure, the Kirillov--Kostant structure. It is characterized
by the property that the Poisson bracket of linear functions
coincides with the Lie bracket of the corresponding elements
of $\g$. Let $G$ be any Lie group whose Lie algebra is $\g$, and
let $T^*G$ be its cotangent bundle, with its canonical
symplectic structure.
Then $\g^*$ may be embedded  as the cotangent space at the
identity, a Lagrangian submanifold of $T^*G$. The Poisson
structure on $\g^*$ is the one that makes the right-invariant
projection $l:T^*G\to \g^*$ a Poisson map. Then $T^*G$
may be canonically quantized: the algebra of
differential operators on $G$ is a quantization of the Poisson algebra
of  functions 
on $T^*G$ and right-invariant differential operators
form a subalgebra which is a quantization of the Poisson
algebra of (polynomial) functions on $\g^*$.

For a general Poisson manifold $M$, the program is to
embed $M$ as a Lagrangian submanifold of a symplectic
manifold $\G$ in such a  way that (deformation, geometric, $\dots)$
quantization of $\G$ descends to a quantization of $M$.
The manifold $\G$ is supposed to be a {\em symplectic groupoid},
an abstraction of the algebraic and geometric properties of
$T^*G$. See Sect.\ \ref{S-4} for the definition of symplectic
groupoids.

The difficulties with this program are, on one side, that
symplectic groupoids do not always exist as smooth manifolds. 
On the other side, it does
not seem to be completely clear in general how to quantize $\G$ 
in such a way that the quantization descends to a quantization
of $M$.

In the meantime, Kontsevich \cite{Ko} found a
different approach to deformation quantization
and constructed  star products for  general Poisson manifolds.

In this paper, we show that (with hindsight)
the program of deformation quantization
based on symplectic groupoids works, albeit in a rather indirect
way.
For each Poisson manifold $M$ we construct a canonical
object $\G$, the phase space of the Poisson sigma model with
target space $M$. The latter is a classical topological field theory.
In its Hamiltonian formulation it is given by an infinite
dimensional Hamiltonian system with constraints. The
constraints generate Hamiltonian vector fields forming an
integrable distribution of tangent  subspaces of
codimension $2\,\mathrm{dim}(M)$ on the constraint surface.
The phase space $\G$ is then the space of leaves of the
corresponding foliation. It carries a natural structure of groupoid,
and also of a symplectic groupoid, in the sense that ``the axioms
would be fulfilled if  $\G$ were a manifold''.   The relation with
the deformation quantization of $M$ is that the same Poisson
sigma model, in its perturbative path integral quantization yields
Kontsevich's deformation quantization formula, as was shown in
\cite{CF1}.

 We study several examples where $\G$ is a manifold   
and also an example, suggested by Weinstein, where it is not: the latter is 
a rotation invariant
Poisson structure on $\R^3$ whose
symplectic leaves are spheres centered at the origin. If the
symplectic area as a function of the radius is not constant
but  has a critical point,
it is known that no symplectic groupoid can 
exist, since it would contradict a theorem of Dazord \cite{D}.
 We show how conical singularities
of $\G$
develop in this case.

In general, the singularities of $\G$ stem from the global structure
of the foliation. However, if we work in the setting of formal
power series, taking a Poisson structure of the form $\epsilon\alpha$
with $\epsilon$ a formal parameter, then a symplectic groupoid
may be constructed \cite{CF2}: it is a formal deformation of the cotangent 
bundle of $M$.

We also note that our $\G$ is related to the 
``local phase space'' of Karasev \cite{Ka}, \cite{KM}.
His construction is based on first order equations
which are essentially our constraint
equation with a special choice of gauge, valid near the
identity elements of the groupoid.

Technically, to work with infinite dimensional manifolds,
we use the framework of manifolds modeled on a Banach(able)
space, for which we
refer to \cite{L}.

The plan of the paper is as follows. In Sect.~\ref{S-2} we introduce
the Poisson sigma model and recall its relation with deformation
quantization. The construction of the
phase space of this model is explained in Sect.~\ref{S-3}.  In
 Sect.~\ref{S-4} we describe the groupoid structure of 
the phase space.

We then turn to examples. In Sect.~\ref{S-5} the basic
examples are presented: in the case of a symplectic manifold $M$
our symplectic groupoid is the fundamental groupoid of $M$ and
in the case of the dual of a Lie algebra it is the cotangent
bundle of the connected, simply connected Lie group
with the given Lie algebra.
In  Sect.~\ref{S-6} we examine a counterexample.

In the last section we study the case of a two-dimensional
Poisson domain and give a construction of a smooth symplectic
groupoid in this case.

\begin{Ack}
We are grateful to L. Tomassini for useful comments, to 
M. Bordemann for interesting discussions and to A. Weinstein for
useful explanations and references to the literature.
\end{Ack}

\section{Poisson sigma model}\lab{S-2}

Let $M$ be a smooth paracompact finite-dimensional
manifold. A Poisson structure on $M$ is a bivector field
$\alpha\in C^\infty(M,\wedge^2TM)$ so that $\{f,g\}=
\alpha(\di f,\di g)$ defines a Lie algebra structure
on the space of smooth functions on $M$. We
assume that a Poisson structure on $M$ is given.
 Let $\Sigma$ be
a two-dimensional oriented compact manifold with boundary.
 We consider an action functional on the space of
vector bundle morphisms $\hat X:T\Sigma\to T^*M$ 
from the tangent bundle of $\Sigma$ to the cotangent bundle
of $M$. Such a map is given by its base map $X:\Sigma\to M$
and a section $\eta$ of $\mathrm{Hom}(T\Sigma,X^*(TM))$:
for $u\in \Sigma, v\in T_u\Sigma$, $\hat X(u,v)=(X(u),\eta(u)v)$.
The action functional is defined on maps obeying the boundary condition
\begin{equation}\label{e-1}
\eta(u)v=0, \qquad u\in\partial \Sigma, \quad v\in T(\partial \Sigma).
\end{equation}
Denote by $\langle\ ,\ \rangle$ the pairing between
the cotangent and tangent space at a point of $M$.
If $X$ is a map from $\Sigma$ to $M$, then this pairing induces
a pairing between the differential forms on $\Sigma$ with values in
the pull-back $X^*(T^*M)$ and the differential forms
on $\Sigma$  with values in $X^* TM$. It
is defined as the pairing of the values and the 
exterior product of differential forms, and takes values in the
differential forms on $\Sigma$.
Then the action functional is 
\[
S(X,\eta)=\int_\Sigma\langle\eta, \di X\rangle
+\frac12\langle\eta,(\alpha\circ X)\eta\rangle.
\]
Here $\eta$, $\di X$ are viewed as one-forms on $\Sigma$
with values in the pull-back of the (co)tangent bundle and
$\alpha(x)$ is viewed as a linear map $T^*_xM\to T_xM$: 
$\alpha(x)=\sum\xi_i\wedge\zeta_i$ is identified with the
map $\beta\mapsto \sum(\xi_i\langle\beta,\zeta_i\rangle
-\zeta_i\langle\beta,\xi_i\rangle)$.
A natural space of vector bundle morphisms in our setting 
 consists of pairs $(X,\eta)$ with
$X$ continuously differentiable and $\eta$ continuous, obeying
the boundary condition \eqref{e-1}. This model was first considered
(in the case of surfaces $\Sigma$ without boundary) in \cite{I},\cite{SS}. 
 
The Feynman path integral for this model with $\Sigma$ a
disk was studied in \cite{CF1}: if $p,q,r$ are three distinct
points on the boundary of the disk, then the semiclassical
expansion of the path integral
\[
f\star g
\,(x)=\int_{X(r)=x}f\left(X(p)\right)
g\left(X(q)\right)e^{\frac i\hbar S(X,\eta)}dXd\eta
\]
around the critical point $X(u)=x$, $\eta=0$ gives 
Kontsevich's star product \cite{Ko} formula.
This action functional is invariant under an infinite dimensional
space of infinitesimal symmetries and the above integral has
to be properly gauge fixed.

Here we want to study the classical part of this story and
formulate the model in the Hamiltonian formalism.

\section{The phase space of the Poisson sigma model}\lab{S-3}

The Hamiltonian formulation of the Poisson sigma model
is obtained by taking $\Sigma$ to be a rectangle $[-T,T]\times I$
with coordinates $(t,u)$ (times and space). The action can
then be written as $S=\int_\Sigma (-\langle\eta_u,\partial_t X\rangle
+\langle \partial_uX +\alpha\eta_u,\eta_t\rangle)
 \di u \di t$. The boundary conditions for $\eta_t$ are
$\eta_t=0$ on $[-T,T]\times \partial I$.
According to the rules of Hamiltonian mechanics,
the first part of this action defines a symplectic structure
on the space of vector
bundle morphisms $TI\to T^*M$ and the coefficient
of the Lagrange multiplier $\eta_t$ is a system of constraints
 that generate a distribution of subspaces spanned by
Hamiltonian vector fields. The phase space of the Poisson sigma
model is then obtained by Hamiltonian reduction, as the set
of integral manifolds of this distribution contained in the set
of zeros of the constraints. It may also be 
expressed as Marsden--Weinstein reduction for the 
symplectic action of an infinite dimensional Lie algebra
 on (a version of)  the cotangent bundle of the path space
$\PM$
of maps $I\to M$.

\subsection{The cotangent bundle of $\PM$}
Let $I$ be the interval $[0,1]$ and $\PM$ the space
of continuously differentiable maps $I\to M$. The tangent
space at $X\in \PM$ is the space of maps $V:I\to TM$
with $V(u)\in T_{X(u)}M$. Let 
$T^*\PM$ be the space of continuous
vector 
bundle morphisms $(X,\eta):TI\to T^*M$ with continuously
differentiable base map
$X:I\to M$. Then $T^*\PM$ is a vector bundle over $\PM$.
The fiber
$T_{X}^*\PM$ at $X$ 
may be thought of as the space of continuous 1-forms on $I$
with values in $X^*(T^*M)$. The vector bundle
 $T^*\PM$ may be thought of as the 
cotangent bundle of $\PM$, since we have a non-degenerate 
pairing  $(\eta,V)\mapsto\int_0^1
\langle\eta(u),V(u)\rangle$ between
$T_X^*\PM$ and $T_X\PM$.
The canonical symplectic form $\omega$
on $T^*\PM$ is defined 
as the differential of the 1-form $\theta_{(X,\eta)}(V)=
-\int_0^1\langle\eta(u),p_* V(u)\rangle$, $V\in T_{(X,\eta)}T^*\PM$,
where $p:T^*\PM\to \PM$ is the bundle projection.

In local coordinates $\hat X$ is described by $n=\mathrm{dim}(M)$ 
functions $X^i\in C^1(I)$ and 
$n$ 1-forms $\eta_i\in C^0(I)du$ on $I$. 
 The symplectic form reads then
\begin{equation}\label{e-omega}
\omega_{\hat X}(\delta_1\hat X,\delta_2\hat X)
=\int_0^{1}
\left( \delta_1X^i\delta_2\eta_i
-
\delta_2X^i\delta_1\eta_i\right).
\end{equation}
We use here and below the Einstein summation convention and do
not write the summation signs for sums over repeated indices.
\subsection{The constraint manifold}
Let $\mathcal{C}$ be the space of solutions of the constraint
equation (``Gauss law'')
\begin{equation}\label{e-Gauss}
\di X(u)+\alpha(X(u))\eta(u)=0,
\end{equation}
with $X$ continuously
differentiable and $\eta$ continuous. This space can be made
into a smooth manifold modeled on the Banach space $\R^n\times
C^0(I,\R^n)$: pick a Riemannian metric on $M$. Then for each $x_0\in M$
and $\eta_0$ a continuous 1-form on $I$ with values
in $T^*_{x_0}M$ sufficiently small, 
there exists a unique solution of \eqref{e-Gauss} 
such that $X(0)=x_0$ and $\eta(u)$ is obtained
from $\eta_0(u)$ by parallel translation for the Levi-Civita connection
along the path $X$. All solutions of \eqref{e-Gauss} may be obtained
this way. Thus, upon choosing local coordinates on a neighborhood
$U\subset M$ of a point and an orthonormal
basis in each tangent space, we have a chart
$\mathcal{C}\supset \mathcal{U}\to\R^n\times C^0(I,\R^n)$.
\subsection{An integrable distribution of subspaces}
Let $\hat X=(X,\eta)$ be a vector bundle morphism
 $TI\to T^*M$ and suppose
$\beta$ is  a continuously differentiable function $I\to T^*M$ such that
$\beta(u)\in T^*_{X(u)}M$, $\forall u\in I$ and $\beta(0)=\beta(1)=0$.
In other words, $\beta$ is in the Banach space 
$C^1_0(I,X^*(T^*M))$ of $C^1$ sections of the pull-back bundle
$X^*(T^*M)$, vanishing at the endpoints. Let
\[
H_\beta=\int_0^1\langle \di X+\alpha\eta,\beta\rangle.
\]
If we vary $\hat X$ in some open set and let $\beta$ depend
on $\hat X$ then $H_\beta$ defines a Hamiltonian vector field
$\xi_\beta$ (``the infinitesimal gauge transformation with gauge parameter
$\beta$'')
on this open set by the rule
\[
\iota_{\xi_\beta}\omega=\di H_\beta.
\]
Here $\iota$ denotes
interior multiplication. This rule makes sense if the dependence
of $\beta$ on $\hat X$ is such that $\di H_\beta$ is in the image
of $\omega$. We show below a way to extend any given 
$\beta\in C^1_0(I,X^*(T^*M))$ in such a way that this holds.
If $\hat X\in\mathcal{C}$,
then $H_\beta$ vanishes and the value of $\di H_\beta$ at $\hat X$
only depends on $\beta$ at $\hat X$. Therefore we have 
for each solution $\hat X$ of \eqref{e-Gauss} a subspace
of the tangent space to the space of vector
 bundle morphisms $TI\to T^*M$
at $\hat X$ spanned by the vectors $\xi_\beta$, $\beta\in
C^1_0(I,X^*(T^*M))$. A formula for $\xi_\beta$ is the following.
Let $\nabla^{TM}$ be a torsion-free connection on $TM$. This connection
induces connections $\nabla^{T^*M}$, $\nabla^{\wedge^2TM}$,
$\nabla^{X^*(T^*M)}$ on the vector bundles $T^*M$, $\wedge^2TM$ over
$M$, and $X^*(T^*M)$ over $I$, respectively. For $x\in M$, 
let $p\in T^*_xM$, $h_{(x,p)}$ denote the
{\em horizontal lift}\/ homomorphism $T^*_xM\to T_{(x,p)}(T^*M)$. It 
maps the tangent vector to a curve $\gamma$ through $x$ to the tangent vector
of the curve $\hat\gamma$ through $(x,p)$ obeying the geodesic equation 
$\nabla^{\gamma^*(T^*M)}\hat\gamma=0$. If $\hat X=(X,\eta)\in T^*\PM$,
then $\xi_\beta(\hat X)$ is the vector bundle morphism $TI\to T(T^*M)$
\begin{eqnarray*}
\xi_\beta(\hat X)(u,v)&=&-h_{(X(u),\eta(u)v)}\left(\alpha(X(u))\beta(u)\right)
                      +\nabla^{X^*(T^*M)}_v\beta(u)\\
 &&              -\langle\beta(u),(\nabla \alpha)(X(u))\eta(u)v\rangle,
\qquad u\in I,\quad v\in T_uI.
\end{eqnarray*}
The last two terms are in $T^*_{X(u)}M$ which is identified
with the vertical tangent space at $(X(u),\eta(u)v)\in T^*M$. 
If $\hat X$ solves \eqref{e-Gauss}, then this expression is independent
of the choice of the connection. It may be more illuminating to write
$\xi_\beta$ in local coordinates: applying $\xi_\beta$ to
the coordinate maps $\hat X\mapsto X^i(u)$, $\hat X\mapsto \eta_i(u)$,
with respect to some choice of coordinates on $M$, gives
\begin{eqnarray}\label{e-gauge}
\xi_\beta X^i(u)&=&
-\alpha^{ij}(X(u))\beta_j(u)\notag
\\
 & & 
\\
\xi_\beta \eta_i(u)&=&\di_u\beta_i(u)+
\partial_i\alpha^{jk}(X(u))\eta_j(u)\beta_k(u)\notag
\end{eqnarray}

\begin{theorem}\label{t-00} Let $\hat X =(X,\eta)\in\mathcal{C}$.
Then the subspace of $T_{\hat X}\PM$ spanned by $\xi_\beta$,
$\beta\in C^1_0(I,X^*(T^*M))$, is a closed subspace of codimension 
$2\,\mathrm{dim}(M)$.
\end{theorem}

\begin{Proof}
For simplicity, we present the proof for $M$ a domain in $\R^n$ and
work with coordinates.
A general tangent vector at a point $(X,\eta)$ of $\mathcal{C}$ 
is a solution
$(\dot X,\dot\eta)$ of the linearization 
\[
\dot X^i(u)+\partial_k\alpha^{ij}\left(X(u)\right)\dot X^k(u)\eta_j(u)
+\alpha^{ij}\left(X(u)\right)\dot\eta_j(u)=0,
\]
of the constraint equation. With our conditions on differentiability,
$(\dot X,\dot\eta)\in C^1(I,\R^n)\oplus C^0(I,\R^n)$,
the map $\beta\to\xi_\beta$ is a continuous linear map from the
Banach space $C^1_0(I,\R^n)$  to $C^1(I,\R^n)\oplus C^0(I,\R^n)$.
It is injective, since $\xi_\beta=0$ implies that $\beta$ obeys
a homogeneous linear first order differential equation with zero
initial condition, and thus vanishes identically. 

Let us describe the image of $\xi$. If $(\dot X,\dot\eta)$ is in the
image then $\dot X(0)=0$ and $\dot \eta$ is of the form
\begin{equation}\label{e-dedb}
\dot\eta(u)=\di_u\beta(u)+A(u)\beta(u),
\end{equation}
for some $\beta\in C^1_0(I,\R^n)$, where $A(u)$ is the matrix
$(\partial_i\alpha^{kj}\eta_j)_{i,j=1\dots,n}$. If $V(u)$ is
the solution of $d_uV(u)=V(u)A(u)$ with $V(0)=1$, then \eqref{e-dedb}
 reads $V(u)\dot\eta(u)=\di_u\left(V(u)\beta(u)\right)$.
Since $\beta(u)$ vanishes at the endpoints, we see that
$\int_IV(u)\dot\eta(u)=0$. Conversely, if $(\dot X,\dot \eta)$ obey
\begin{equation}\label{e-qq}
\dot X(0)=0,\qquad \int_IV(u)\dot\eta(u)=0,
\end{equation}
then $(\dot X,\dot\eta)=\xi_\beta$, with $\beta(u)=V(u)^{-1}\int_0^uV(u')
\eta(u')$.

The image is thus described as the common kernel 
\eqref{e-qq} of $2n$ linearly
independent continuous linear functions, and is thus closed of
codimension $2n$.
\end{Proof}

The next step is to show that the distribution of subspaces in the
tangent bundle to the space of solutions of $\eqref{e-Gauss}$ is integrable
and thus defines a  foliation of codimension $2\,\mathrm{dim}(M)$.
This is best seen by interpreting the leaves as orbits of 
a gauge group which we introduce in the next section.

\subsection{The Lie algebra and its action on the cotangent bundle}
The Lie algebra acting on $T^*\PM$ is obtained from the
Lie algebra of 1-forms $\Omega^1(M)$ with the Koszul Lie bracket. 
This bracket is defined by
\[
 [\beta,\gamma]=
\di \langle\beta,\alpha\gamma\rangle
-\iota_{\alpha\beta} \di \gamma
+\iota_{\alpha\gamma} \di \beta,
\]
for any $\beta,\gamma\in\Omega^1(M)$. 
 In local coordinates
$\alpha=\frac12
\alpha^{ij}\frac\partial{\partial x^i}\wedge
\frac\partial{\partial x^j}$, 
$\beta=beta_i\di x^i$, 
$\gamma=\gamma_i\di x^i$, with $\partial_i=\partial/\partial x^i$,
\[
[\beta,\gamma]=
\left(
\partial_i\alpha^{jk}\beta_j\gamma_k
+
\alpha^{jk}\partial_j\beta_i\gamma_k
+
\alpha^{jk}\beta_j\partial_k\gamma_i
\right)\di x^i.
\]
This bracket obeys the Jacobi identity
as a consequence of the Jacobi identity for $\alpha$.
Let $P_0\Omega^1(M)$ be the Lie algebra of continuously differentiable
maps $I\to \Omega^1(M)$
such that $\beta(0)=\beta(1)=0$,
 with bracket $[\beta,\gamma](u)
=[\beta(u),\gamma(u)]$.

If $\beta\in P_0\Omega^1(M)$, let 
\[
H_\beta(X,\eta)=\int_I\langle \di X(u)+\alpha(X(u))\eta(u),
\beta(X(u),u)\rangle.
\]
Recall that if $H$ is a smooth function on a symplectic manifold, then
a vector field $\xi$ is called Hamiltonian vector field generated by $H$
if $\iota_\xi\omega=\di H$. Such a vector field, it it exists, is unique.
In the infinite dimensional setting existence is not guaranteed in
general. 
\begin{theorem}\label{t-1}
\ 
\begin{enumerate}
\item[(i)] For each $\beta\in P_0\Omega^1(M)$
there exists a Hamiltonian vector field $\xi_\beta$ generated
by $H_\beta$.
\item[(ii)] The Lie algebra $P_0\Omega^1(M)$ 
acts on $T^*\PM$ by the Hamiltonian
vector fields $\xi_\beta$, i.e., $\beta\mapsto \xi_\beta$ is
a Lie algebra homomorphism.
\item[(iii)] The map $\mu:T^*\PM
\to P_0\Omega^1(M)^*$ with $\langle\mu(X,\eta),\beta\rangle
=H_\beta(X,\eta)$ is an equivariant moment map for this action.
\end{enumerate}
\end{theorem}
\begin{Proof}
(i) By using a partition of unity, we may restrict ourselves to
$\beta$ with support in a coordinate neighborhood of $M$, and
use local coordinates. If $\dot X^i(u),\dot\eta_i(u)$ are the coordinates
of a vector field $\zeta$ on $T^*\PM$ then, with the
abbreviations $\beta_i=\beta_i(X(u),u)$, $\alpha^{ij}=\alpha^{ij}(X(u))$,
\begin{eqnarray*}
\di H_\beta(\zeta)&=&
\int_0^1\left((\di_u\dot X^i+\partial_k\alpha^{ij}\dot X^k
\eta_j+\alpha^{ij}
\dot\eta_j)
             \beta_i
+C^j\partial_i\beta_j\dot X^i\right)
\\
&=&
\int_0^1\dot X^i\left(-\di_u\beta_i
+\partial_i\alpha^{jk}\eta_k\beta_j+C^j\partial_i\beta_j\right)
+\int_0^1\dot \eta_i\alpha^{ji}
             \beta_j.
\end{eqnarray*}
The term with $C^j=\di_u{X^i}+\alpha^{ij}\eta_j$ vanishes on $\mathcal{C}$.
Here $\di_u$ is the (total) differential with respect to the
coordinate $u$ on the interval.
We may then read off the coordinates $\delta_\beta X^i,\delta_\beta\eta_i$
of $\xi_\beta$, and at the same time
show that they exist, from the defining relation
$\omega(\xi_\beta,\zeta)=\di H_\beta(\zeta)$, where $\omega$ is
given by \eqref{e-omega}. We obtain
\begin{eqnarray}\label{e-mezzogiorno}
\delta_\beta X^i(u)&=&-\alpha^{ij}(X(u))\beta_j(X(u),u)\notag
\\
\delta_\beta \eta_i(u)&=&\di_u\beta_i(X(u),u)+\partial_i
\alpha^{jk}(X(u))\eta_j(u)\beta_k(X(u),u)
\\
&&-C^j(u)\partial_i\beta_j(X(u),u).\notag
\end{eqnarray}

\noindent(ii) is a consequence of (iii)

\noindent(iii) The statement amounts to the identity
$H_{[\beta,\gamma]}=\xi_\beta H_\gamma$, which we may
again check in local coordinates.
We have
\begin{eqnarray*}
\xi_\beta H_\gamma(X,\eta)&=&\xi_\beta\int_0^1C^i\gamma_i
\\
&=&-\int_0^1\biggl(C^i\partial_k\gamma_i\alpha^{kl}\beta_l
-\partial_k\alpha^{ij}\di_u X^k{}\beta_j\gamma_i
\\
&&
-\partial_k\alpha^{ij}\alpha^{kl}\beta_k\eta_j\gamma_i
+\alpha^{ij}\partial_j\alpha^{kl}\eta_j\beta_l-\alpha^{ij}C^k\partial_j\beta_k
\gamma_i\biggr).
\end{eqnarray*}
By combining terms with the Jacobi identity, we arrive at the
formula
\[
\xi_\beta  H_\gamma
=\int_0^1C^i(\partial_i\alpha^{jk}\beta_j\gamma_k
+
\alpha^{jk}\partial_j\beta_i\gamma_k
+
\alpha^{jk}\beta_j\partial_k\gamma_i
)=H_{[\beta,\gamma]}.
\]
\end{Proof}

\subsection{The phase space} The set $\mu^{-1}(0)$ of zeros of
the moment map is  the constraint manifold $\mathcal{C}$.
One would like to define the phase space as the Marsden--Weinstein
symplectic quotient $T^*\PM\dslash H=\mathcal{C}/H$.
The gauge group $H$  is the group of symplectic
diffeomorphisms generated
by the flows of the Hamiltonian vector fields $\xi_\beta$.
The trouble is that not only the manifold is infinite dimensional,
but the action of the group is far from being nice, and one should
not expect to have a good quotient.

However, locally the orbits form a smooth foliation:

\begin{theorem}
The distribution of tangent subspaces of $\mathcal{C}$ spanned by 
$\xi_\beta$, $\beta\in C^1_0(I,X^*(T^*M))$ is integrable. Its
integral manifolds are smooth of codimension $2\,\mathrm{dim}(M)$ and
are the orbits of $H$.
\end{theorem}

\begin{Proof}
We present the proof in the case where $M$ is a domain in $\R^n$.
The general case is treated in a similar but more cumbersome way.

Let $V_{(X,\eta)}$ be the subspace of $T_{(X,\eta)}\mathcal{C}$.
 spanned by $\xi_\beta$,
$\beta\in C^1_0(I,X^*(T^*M))$. These vector spaces form a
smooth subbundle of the tangent bundle: locally over a 
neighborhood $\mathcal{U}\subset\mathcal{C}$, this subbundle
is the image of the smooth vector bundle morphism
\begin{eqnarray*}
\mathcal{U}\times C_0^1(I,\R^n)&\to& T\mathcal{U} \\
\left((X,\eta),\beta\right)&\mapsto& \left((X,\eta),\xi_\beta(X,\eta)\right).
\end{eqnarray*}
By Theorem \ref{t-00}, in each fiber this is an injective
linear continuous map with closed image of codimension $2n$.

Now the integrability follows from the Frobenius theorem (see
\cite{L}, Chapter VI, for a proof valid in the infinite dimensional
setting): every $\beta\in
C^1_0(I,X^*(T^*M))$ may be extended to an element of $P_0\Omega^1(M)$:
in coordinates, take $\beta_i(x,u)$ independent of $x$. It then
follows from Theorem \ref{t-1} (ii), that $[\xi_\beta,\xi_\gamma]=
\xi_{[\beta,\gamma]}$, which implies the Frobenius integrability 
criterion.

The
fact that the integral manifolds are orbits of $H$, follows from the
fact that $V_{(X,\eta)}$ coincides with the space spanned by
the restriction to $(X,\eta)$ of Hamiltonian vector fields
generated by $H_\beta$, $\beta\in P_0\Omega^1(M)$. 
\end{Proof}

\section{The symplectic groupoid structure on $T^*\PM\dslash H$}
\lab{S-4}
A symplectic manifold $\G$ with
symplectic form $\omega_\G$ is called {\em symplectic groupoid}
for a Poisson manifold $M$ if we have an injection
 $j:M\hookrightarrow
\G$, two surjections $l,r:\G\to M$, a composition law
$g,h\mapsto g{\sbullet} h$ defined if $g,h\in\G$ and $r(g)=l(h)$
obeying a set of axioms. The first axioms say that $\G$ is a
groupoid, i.e., denoting $\G_{x,y}=l^{-1}(x)\cap r^{-1}(y)$,
\begin{enumerate}
\item[(i)] $l\circ j=r\circ j=\mathrm{id}_M$.
\item[(ii)] If $g\in \G_{x,y}$ and $h\in\G_{y,z}$, then
$g\sbullet h\in \G_{x,z}$.
\item[(iii)] $j(x)\sbullet g=g\sbullet j(y)=g$, if $g\in\G_{x,y}$.
\item[(iv)] To each $g\in\G_{x,y}$ there exists an inverse $g^{-1}\in\G_{y,x}$
such that $g\sbullet g^{-1}=j(x)$.
\item[(v)]  The composition law is associative:
$(g\sbullet h)\sbullet k=g\sbullet(h\sbullet k)$ whenever
defined.
\end{enumerate}
In the language of categories, these axioms say that $\G$ is
the set of morphisms of a category in which all morphisms
are isomorphisms. $M$ is the set of objects and $j(M)$ the
set of identity morphisms.
It follows from the axioms that $g^{-1}$ is uniquely determined
by $g$ and that
$g^{-1}\sbullet g=j(y)$, if $g\in\G_{x,y}$. 

The next axioms
relate to the symplectic and Poisson structure. A smooth  map 
$\phi:M_1\to M_2$ between Poisson manifolds
is called Poisson if $\{f,g\}\circ\phi=\{f\circ \phi,g\circ\phi\}$, for
all $f,g\in C^\infty(M_2)$. It is called anti-Poisson if 
 $\{f,g\}\circ\phi=-\{f\circ \phi,g\circ\phi\}$.
Then the remaining axioms are:
\begin{enumerate}
\item[(vi)] $j$ is a smooth embedding, $l,r$ are smooth submersions,
the composition and inverse maps are smooth.
\item[(vii)] $j(M)$ is a Lagrangian submanifold. In particular
$\mathrm{dim}(\G)=2\,\mathrm{dim}(M)$.
\item[(viii)] $l$ is a Poisson map and $r$ is an anti-Poisson map.
\item[(ix)] Let $P:\G_0\subset\G\times \G\to \G$ be the
composition law on $\G_0=\{(g,h)\in\G\,|\,r(g)=l(h)\}$, and
$\pi_1,\phi_2:\G\times\G\to\G$ denote the projections onto
the first and second factor. Then $P^*\omega_\G=\pi_1^*\omega_\G
+\pi_2^*\omega_\G$.
\item[(x)] $g\mapsto g^{-1}$ is an anti-Poisson map.
\end{enumerate}
The basic example is the following:
\begin{Exa}\label{Liesg}
Let $M=\mathfrak{g}^*$ be the dual space to a Lie algebra $\mathfrak g$
with Kirillov--Kostant Poisson structure. For any Lie group
$G$ with Lie algebra isomorphic to $\g$,
 we have the inclusion $j:\g^*\to T^*G$ of
$\g^*$ as the cotangent space at the identity  $e$ and
projections $l,r:T^*G\to\g^*$ sending the cotangent
space at each point to the cotangent space at the identity by
left (right) translation. If $L_g,R_g:G\to G$, with
$L_g(h)=gh$, $R_g(h)=hg$, denote
the left and right translation by $g$, we have $l(g,a)=\di R_g(e)^*a$,
$r(g,a)=\di L_g(e)^*b$, ($g\in G$, $a\in T^*_{g}G$). 
If $r(g,a)=l(h,b)$, the composition law is $(g,a)\sbullet
(h,b)=(gh,c)$ with $c=(\di R_h(g)^*)^{-1}a=(\di L_g(h)^*)^{-1}b$.
 
A more explicit description is obtained by identifying
$T^*G$ with $\g^*\times G$ via $(g,a)\mapsto (\di R_g(e)^*a,g)$,
see \ref{S-dLa} below.
\end{Exa}

\subsection{The groupoid structure}
The algebraic 
groupoid structure of $\G=\mathcal{C}/H$ can be naturally
defined in terms of composition of paths. 
We have an inclusion $j:M\hookrightarrow \G$ sending
a point $x$ to the class of the constant solution $X(u)=x,
\eta(u)=0$.  Let 
$l,r:T^*\PM\to M$ be the maps
\begin{equation}\label{lrXeta}
l(X,\eta)=X(0),\quad r(X,\eta)=X(1).
\end{equation}
These maps are $H$-invariant, since the symmetries preserve
the endpoints, hence they descend to maps $l,r:\G\to M$,
and it is clear that axiom (i) holds. 

\subsection{Composition law and inverses}
To define the composition
law we need to do some adjustments at the endpoints:
Let $H_0$ be the subgroup of $H$ generated by the flows of
the vector fields $\xi_\beta$ such that $\di \beta(0)=\di \beta(1)=0$.

\begin{lemma}\label{l-1}
In each equivalence class $[(X,\eta)]$ in $\G=\mathcal{C}/H$ there
exists a representative with $\eta(0)=\eta(1)=0$. Any two representatives
with this property can be related by an element of $H_0$.
\end{lemma}

\begin{Proof}
Let $(X,\eta)\in \mathcal{C}$. To obtain a representative with
$\eta(0)=0$ we perform a gauge transformation obtained as the
flow of a vector field $\xi_\beta$ with $\beta$ supported on
a small neighborhood $I_0$ of $0\in I$. Small means here that 
$X(u)$ lies in a coordinate neighborhood $U$ of $M$ for $u\in I_0$.
Then the gauge transformation may be described in local coordinates.
The problem is then to find continuously differentiable
functions $\beta_i(u)$ supported on
$I_0$ with $\beta_i(0)=0$, so that the solution to the system
\begin{eqnarray*}
\frac{\partial}{\partial s}X^i(u,s)&=&- \alpha^{ij}\left(X(u,s)\right)
\beta_j(u),\\
\frac{\partial}{\partial s}\eta_i(u,s)&=&\di_u\beta_i(u)
+ \partial_i\alpha^{jk}\left(X(u,s)\right)
\eta_j(u,s)\beta_k(u),
\end{eqnarray*}
with initial condition $(X,\eta)$ at $s=0$ (a) exists with $X$ in
 $U$
for all $s\in[0,1]$, and (b) obeys $\eta_i(0,1)=0$. 
A sufficient condition for (a) is that $|\beta_i(u)|<\delta$ for
some $\delta>0$ and all $u\in I_0$: if this bound holds with
sufficiently small $\delta$, then
the first equation has a solution which is close to $X$ and thus
remains in $U$. Given $X(u,s)$, the second equation is linear
for $\eta$ and thus has a solution for all $s\in[0,1]$. 
To achieve (b), let $a_i=\eta_i(0)$
and choose $\beta_i(u)$ so that $\beta_i(u)=-a_iu+O(u^2)$. Then
$\eta_i(0,s)=(1-s)a_i$ vanishes for $s=1$. The same procedure
may be applied at the other end of $I$.

Suppose now that $\hat X^{(0)}=(X^{(0)},\eta^{(0)})$ and $\hat X^{(1)}=(X^{(1)},\eta^{(1)})$
are two representatives of a class in $\G$, obeying the
condition $\eta^{(0)}(u)=\eta^{(1)}(u)=0$ for $u=0,1$. These representatives
are related by an element of $H$, which is a product of a
finite number  $k$ of flows
of vector fields of the form $\xi_\gamma$, $\gamma\in P_0\Omega^1(M)$.
Let us first assume that $k=1$. Then we have a smooth path
$s\mapsto\hat X_s$ in $\mathcal{C}$, such that $\hat X_{s=0}=\hat X^{(0)}$,
$\hat X_{s=1}=\hat X^{(1)}$, and $\di \hat X_s/\di s=\xi_{\gamma}(\hat X_s)$.
We now repeat the procedure of the first part of the proof, for each
$s\in[0,1]$. Let $x_0=\hat X_s(0)$, which is independent of $s$, and
choose coordinates in a neighborhood $U\subset M$ of $x_0$. 
Let $\beta_s\in P_0\Omega^1(U)$ be such that (a)
$\beta_{s,i}(x,u)=-a_{s,i}u+O(u^2)$
 $(u\to 0$), where $\eta_{s,i}(0)=a_{s,i}\di u$;
(b) $\beta_s(x,u)=0$ if $u\geq\delta'$, for some sufficiently small
$\delta'>0$, (c)$\max|\beta_{s,i}|$ is sufficiently small.
Then the flow of $\xi_{\beta_s}$ exists on $\mathcal{ U}=\{
\hat X\in\mathcal{U}\,|\,\hat X(u)\in U,$ for $u\le\delta'\}$
for all times in $[0,1]$. 
Applying this flow to $\hat X_s$ we obtain
a two-parameter family $\hat X_{s,\sigma}$, $(s,\sigma)\in[0,1]^2$, in
$\mathcal{C}$,
differing from $\hat X_s$ only in some small neighborhood of
$0\in I$,
 such that
$\hat X_{s,0}=\hat X_s$ and $\partial\hat X_{s,\sigma}/\partial\sigma=
\xi_{\beta_s}(\hat X_{s,\sigma})$. By construction, we have
\begin{enumerate}
\item[(i)] $\hat X_{0,\sigma}=\hat X^{(0)}$, $\hat X_{1,\sigma}=\hat X^{(1)}$,
for all $\sigma\in[0,1]$.\label{item1}
\item[(ii)] $\hat X_{s,1}=(X_{s,1},\eta_{s,1})$ with $\eta_{s,1}(0)=0$.
\end{enumerate}\label{item2}
Since $\xi_{\beta_s}$, $\xi_{\gamma}$ belong to an integrable distribution
of tangent subspaces, any curve $t\mapsto\hat X_{s(t),\sigma(t)}$
is in the integral manifold passing through $\hat X^{(0)}$. In
particular, $s\mapsto\hat X_{s,1}$ defines by (i), (ii)
a curve of points related by a transformation of $H$. Since
$\eta_{s,1}(0)=0$, this transformation is in $H_0$. The
same argument applies to the other endpoint, and for $k\ge1$ one
applies this construction $k$ times.
\end{Proof}

Then we may define the composition law
$[(X,\eta)]=[(X_1,\eta_1)]\sbullet[(X_2,\eta_2)]$  in $\G$
by choosing representatives as in Lemma \ref{l-1} and setting
\begin{equation}\label{compXeta}
\begin{split}
X(u)&=\left\{
\begin{array}{rl}
X_1(2u),& 0\leq u\leq\frac12,\\
X_2(2u-1),& \frac12\leq u\leq1,
\end{array}
\right.
\\
\eta(u)&=\left\{
\begin{array}{rl}
2\,\eta_{1u}(2u)\,\di u,& 0\leq u\leq\frac12,\\
2\,\eta_{2u}(2u-1)\,\di u,& \frac12\leq u\leq1.
\end{array}
\right.
\end{split}
\end{equation}
(we write $\eta_i(u)=\eta_{iu}(u)\,\di u$, $i=1,2$), provided 
$X_1(1)=X_2(0)$.

Lemma \ref{l-1} ensures that $\eta$ is continuous. $X$ is continuously
differentiable since the derivatives of $X_1$, $X_2$  at the endpoints match
--- they vanish,
as a consequence of \eqref{e-Gauss}. It is immediate to check that $(X,\eta)$
obeys the constraint equation if $(X_1,\eta_1)$, $(X_2,\eta_2)$ do.
Therefore the composition is well-defined at the level of representatives.
By the second part of Lemma \ref{l-1}, the class of $(X,\eta)$ is independent
of the choice of representatives: infinitesimal transformations of
$(X_1,\eta_1)$, $(X_2,\eta_2)$ associated to $\beta_1,\beta_2:I\to T^*M$
and obeying $\di \beta_1(1)=\di \beta_2(0)$ match at the end points to give
to give an infinitesimal transformation of $(X,\eta)$ associated to
\[
\beta(u)=\left\{
\begin{array}{rl}
\beta_1(2u),& 0\leq u\leq\frac12,\\
\beta_2(2u-1),& \frac12\leq u\leq1,
\end{array}
\right.
\]
which is a differentiable function $\beta:I\to T^*M$ 
with $\beta(u)\in T^*_{X(u)}M$.

For $u\in[0,1]$, let
$\theta(u)=1-u$. If $\hat X=(X,\eta)$  obeys the constraint equation 
\eqref{e-Gauss},
then $\theta^*\hat X=(X\circ\theta,\theta^*\eta)$ also does. Moreover the endpoints
of the path $X$ are interchanged under this map. If $\beta$ is a section
of $X^*(T^*M)$ then $\beta\circ\theta$ is a section of $(X\circ\theta)^*(T^*M)$,
and $\theta^*(\xi_{\beta}\hat X)=\xi_{\beta\circ\theta}\theta^*\hat X$.
Therefore 
\[
[(X,\eta)]\mapsto [(X,\eta)]^{-1}=[(X\circ\theta,\theta^*\eta)]
\]
is a well-defined map from $\G$ to $\G$.

\begin{theorem}
 $\G$  obeys
axioms (i)--(v)
\end{theorem}

The idea of the proof of the associativity is based on the fact
that the composition law is associative up to reparametrization
of $I$. But
it turns out that reparametrizations are special gauge transformations:
indeed, if an infinitesimal reparametrization is given by a vector
field $\epsilon$ on $I$ vanishing at the endpoints, then the variation
of a solution $\hat X$ of the constraint equation is an infinitesimal
transformation with parameter
$\beta(u)=\eta(u)\epsilon(u)$, provided $\eta$ is differentiable.
 Similarly, to prove that $g\sbullet g^{-1}=j(x)$
one uses the fact that $g\sbullet g^{-1}$ is the class
of a point $(X,\eta)\in\mathcal{C}$ so that $(X\circ\theta,\theta^*\eta)=(X,\eta)$.
Let
\[
\beta(u)=\beta_{\hat X}(u)=\left\{
\begin{array}{rl}
u\,\eta_u(u),& 0\leq u\leq\frac12,\\
(1-u)\,\eta_u(u),& \frac12\leq u\leq1.
\end{array}
\right.
\]
Here again, we have to assume that $\eta$ is differentiable.
Then the solution of $\frac {\di \hat X_s}{\di s}=\xi_{\beta} $ 
with initial condition $\hat X_0=\hat X$ is
a solution $\hat X_s$ of \eqref{e-Gauss} which is also symmetric
with respect to $\theta$. For $s=1$, $X_1(u)=x$ is constant. Thus
$\eta$ is in the kernel of $\alpha(x)$ and obeys $\theta^*\eta=\eta$.
But Ker$\,\alpha(x)$ is naturally a Lie algebra (the Lie bracket
$[\beta,\gamma]$ is the Koszul bracket for any two
1-forms coinciding with $\beta$, $\gamma$ at $x$). 
Then $\di+\eta$ has the interpretation of a connection
on a trivial vector bundle over $I$. 
Infinitesimal
gauge transformations preserving the condition $X(u)=x$ are
infinitesimal gauge transformations in the usual gauge theory
sense. In particular a connection with $\theta^*\eta=\eta$
is gauge equivalent (with a gauge transformation which
is trivial at the endpoints) to the trivial connection.

Technically, these operations are possible thanks to the
\begin{lemma}
In every class $[(X,\eta)]\in\G$ there exists a representative
so that $X$ and $\eta$ are smooth maps.
\end{lemma}

\begin{Proof}
Let $(X,\eta)\in \mathcal{C}$. Let  us divide the interval $I$
into an odd number
$k\geq3$ of parts $I_j=[j/k,(j+1)/k]$, $0\leq j\leq k-1$, in such
a way that $X(I_j)$ is contained in a coordinate neighborhood
of $M$. To find a smooth representative, we perform a
sequence of gauge transformations. Each of these gauge
transformation is generated by parameters $\beta$ with
support in a small  neighborhood of an interval $I_j$. 
Such gauge transformations only affect $(X,\eta)$ in such
a neighborhood. Therefore we may describe them in
local coordinates by the formula \eqref{e-gauge}.

Let $\etasm\in C^\infty(I,X^*(T^*M))\otimes\Omega^1(I)$
be a smooth section, $C^0$-close to $\eta$.

As a first step, we show that $\eta$ may be taken to be equal
to $\etasm$ 
on $I_0$. Let  for $s\in[0,1]$,
 $\eta_i(u,s)=s\,\etasm_i(u)+(1-s)\,\eta_i(u)$, 
$1\leq i\leq \mathrm{dim}(M)$,
$u\in I_0$. 
Let $X(u,s)$ be the solution of the constraint equation on $I_0$
\[
\frac\partial{\partial u}{X^i}(u,s)+\alpha^{ij}(X(u,s))\eta_j(u,s)=0,
\]
 with $X^i(0,s)=X^i(0)$. This equations has a unique solution on $I_0$
if $\etasm$ is sufficiently close to $\eta$.
Let $\beta$ be the solution of the linear differential equation
\begin{equation}\label{e-mezzanotte}
\frac\partial{\partial u}\beta_i(u,s)+
\partial_i\alpha^{jk}(X(u,s))\eta_j(u,s)\beta_k(u,s)=
\frac\partial{\partial s}\eta_i(u,s),
\end{equation}
on $I_0$ with initial condition $\beta_i(0,s)=0$. Extend $\beta_i(u,s)$
to a function on $I$ vanishing outside some small neighborhood
of $I_0$. Then $\beta_i(\cdot,s)$ is the local coordinate expression of 
a section $\beta_s\in C^1_0(I,X^*(T^*M))$ with support in a 
neighborhood of $I_0$. It may be taken to be small in the
$C^1$-topology if $\etasm$ is close to $\eta$ in the $C^0$-topology. The flow of the vector field $\xi_{\beta_s}$,
$0\leq s\leq1 $, is then a gauge transformation that transforms $(X,\eta)$
into a solution $(\tilde X,\tilde \eta)$ coinciding with $(X,\eta)$
outside a neighborhood of $I_0$ and such that $\tilde\eta=\etasm$ on
$I_0$.

This step may be repeated on $I_1$, $I_2$, \dots, until we
get a representative $(\tilde X,\tilde\eta)$ with $\tilde\eta=\etasm$
on $[0,(k+1)/2k]$. We then repeat the same step integrating
\eqref{e-mezzanotte} backwards,
starting from the last interval $I_{k-1}$, and continuing with
$I_{k-2}, \dots$, until we arrive at the middle interval $I_{(k-1)/2}$.
At this point the representative $(\tilde X,\tilde \eta)$ has
$\tilde\eta=\etasm$ except on some small interval $I'$ in the middle of $I$. 
We apply once more our step to a slightly bigger interval $I''$
including
$I'$. Then the solution $\beta$ of \eqref{e-mezzanotte} is a
smooth function of $u$ on $I''\setminus I'$ and may be extended
to a section in $C^1_0(I,X^*(T^*M))$ which is smooth outside $I'$.
The resulting representative $(\bar X,\bar\eta)$ of the
class $[(X,\eta)]$ has then $\bar\eta$ smooth. Then also
$\bar X$, as a solution of \eqref{e-Gauss}, is smooth.
\end{Proof}

\subsection{Symplectic structure}
To formulate axioms (vi)--(x) we need $\G$ to be a manifold, which
is not always the case, as we shall see below. 

So we assume that $\G$ is a manifold, or more precisely:

\begin{Ass}\label{A-1} There exists a smooth manifold $\G$ and a
smooth submersion $\pi:\mathcal{C}\to \G$ whose fibers are
the $H$-orbits.
\end{Ass}

Below we give examples where this assumption holds and examples
where it does not.

The symplectic structure $\omega_\G$ on $\G$ is constructed in the
usual way: the point is that the symplectic 2-form $\omega$ on $T^*\PM$
restricts to an $H$-invariant
closed 2-form on $\mathcal{C}$ whose null spaces are the tangent spaces to the
orbits. This implies that
\[
\omega_\G(x)(\xi,\zeta)=\omega(\hat X)(\hat\xi,\hat\zeta),
\qquad \xi,\zeta\in T_x\G,
\]
is independent of the choice of $\hat X\in\mathcal{C}$ such that $\pi(\hat X)=x$
or of $\hat\xi,\hat\zeta\in T_{\hat X}\mathcal{C}$ projecting to $\xi,\zeta$,
and defines a symplectic 2-form on $\G$.

\begin{theorem}
Under Assumption \ref{A-1},
$\G$ is a symplectic groupoid
for $M$.
\end{theorem} 

\begin{Proof}
The non-trivial assertions are (viii)-(x). Let us prove that
the left projection $l$ is a Poisson map. Let $\mathcal{U}$
be some small neighborhood in $T^*\PM$ of a point $\hat X_0\in\mathcal{C}$,
and choose local coordinates on $M$ around $X_0(0)$. Then it is
sufficient to show that the coordinates $l^i$ of $l$ obey 
$\{l^i,l^j\}_\G=\alpha^{ij}\circ l$. Let $\psi(u)\,\di u$ 
be any smooth 1-form
on $I$ with support in a small neighborhood of $0$ and such that
$\int_0^1\psi(u)\di u=1$. Then the function $\mathcal{U}\to \R$
\[
l^i_\psi:\hat X\mapsto\int_0^1\left( X^i(u)+\int_0^u\alpha^{ij}\left(
X(v)\right)\eta_j(v)\right)\psi(u)\, \di u
\]
restrict to $l^i$ on $\mathcal{U}\cap\mathcal{C}$ 
(with the support condition on $\psi$,
this local coordinate expression  makes sense).
The main property of this extension of $l^i$ is that its differential
lies in the image of $\omega$ and thus generates a local Hamiltonian
vector field $\xi^i$. Therefore we may compute $\{l^i,l^j\}_\G$ as
the Poisson bracket $\{l^i_\psi,l^j_\psi\}$ on $T^*\PM$, which
is then $\xi^il^j$. The differential of $l^i_\psi$ applied to
a vector field $\zeta$ with coordinates $\dot X^j,\dot\eta_j$ is
\begin{eqnarray*}
\di l^i_\psi(\zeta)&=&
\int_0^1\dot X^j(u)\left(
\delta_{ij}\di u+\partial_j\alpha^{ik}\left(X(u)\right)\eta_k(u)
\int_u^1\psi(v)\,\di v\right)\\
&&+\int_0^1\dot\eta_j(u)\alpha^{ij}\left(X(u)\right)\int_u^1\psi(v)\,\di v.
\end{eqnarray*}
The vector field $\xi^i$, solution of $\omega(\xi^i,\zeta)=\di l^i(\zeta)$
has then coordinates $\delta^iX^j,\delta^i\eta_j$ with
\[
\delta^i X^j(u)=\alpha^{ij}\left(X(u)\right)\int_u^1\psi(v).
\]
We do not need $\delta^i\eta_j$. Then
\[
\{l^i,l^j\}_\G([\hat X])=\xi^il^j(\hat X)=\delta^iX^j(0)=\alpha^{ij}(l(\hat X)),
\]
as was to be shown.
To prove (ix) we notice that the integral defining $\omega$ at the
product of two solutions is the sum of the integrals defining 
$\omega$ at the two solutions.
Axiom (x) follows from the fact that the inversion changes the
sign of the symplectic form, as can be seen by changing $u$ to
$1-u$ in the integral defining $\omega$.
\end{Proof}


\section{Basic examples}\lab{S-5}
In this and in the next sections we discuss some examples. 
To fix the notations,
we will always denote by $u$ the variable in our space interval $I=[0,1]$.
When considering a flow generated by symmetries, we will denote the flow
parameter by $s$. Finally, we will use a prime to indicate derivatives w.r.t.\ 
$u$ and a dot for derivatives w.r.t.\ $s$.

\subsection{Trivial Poisson structures}
Let us consider a manifold $M$ with Poisson bivector field
$\alpha=0$.

In this case, the ``Gauss law'' selects the constant maps $X:I\to M$.

Let $X(u)=\xi\in M$ be such a solution. 
Then the corresponding bundle map $\hat X$
is given by $X$ and a continuous one-form $\eta$ on $I$ that takes value in 
$\T^*_\xi M$.

The infinitesimal symmetries are given by
\[
\delta\eta=\dd\beta,
\]
with $\beta:I\to\T^*_\xi M$, $\beta|_{\de I}=0$.

If we define $\pi:=\int_I \eta \in\T^*_\xi M$, then
for a given solution we have the well-defined map $i:\calG\to\T^*M$,
which maps $\hat X$ into $(X(0),\pi)$.

We can invert this mapping by defining $j:\T^*M\to\calG$ as follows:
$j(g)$, $g\in\T^*M$, is the constant morphism $\hat X(u)=g,\,\forall u\in I$.

An immediate check shows that $i\circ j=\mathrm{id}$. 

We can also prove that $j\circ i=\mathrm{id}$. In fact, let $\hat X\in\calG$. 
Then ${\Hat{\Tilde X}}:=j\circ i(\hat X)$ is a solution with
$\Tilde X = X$, and $\int_I\tilde\eta=\int_I\eta$.
Denoting by $I^u$ the path in $I$ from the lower boundary till
a point $u$, we can then define
$\beta(u)=\int_{I^u}(\tilde\eta-\eta)$, which
is an allowed symmetry generator.

Next we consider the following path of $\T_X^*M$-valued one-forms
\[
\eta_s:= s\,\tilde\eta + (1-s)\,\eta,
\qquad s\in[0,1].
\]
Finally, we have 
\[
\dot\eta_s=\dd\beta,
\]
so that ${\Hat{\Tilde X}}$ is equivalent to $\hat X$.
We have then proved the following
\begin{Thm}
The phase space $\calG$ for a trivial Poisson structure on $M$ is
diffeomorphic to $\T^*M$.
\end{Thm}
It is well-known that $\T^*M$ is a symplectic groupoid for $M$.
The two projections $l$ and $r$ coincide with the natural projection
$T^*M\to M$, while the product is given by
\[
(\xi,\pi)\sbullet(\xi,\pi')=(\xi,\pi+\pi').
\]

\subsection{The symplectic case}
Since now the Poisson bivector field is nondegenerate, the Gauss
law allows to completely determine the bundle morphism 
$\Hat X\colon \T I\to \T^*M$ in terms of the base map $X$:
\[
\Hat X = -\alpha^{-1}(\dd X).
\]

As for $X$, the infinitesimal symmetries are now all infinitesimal
diffeomorphisms of the target that fix the endpoints of $X(I)$.

Thus, the space of solutions modulo symmetries coincides with the 
fundamental groupoid of $M$.

In the case when $M$ is simply connected we can further identify
$\calG$ with $M\times\Bar M$, where $\Bar M$ denotes $M$ with opposite
symplectic structure. The product is then simply
\[
(x,y)\sbullet(y,z)=(y,z).
\]

In the general case, a point in $\calG$ is given by a pair of points 
$x$ and $y$ in $M$ together with a class $c$ of homotopic paths with
fixed endpoints in $x$ and $y$. The product is then
\[
(x,y,c)\sbullet(y,z,c')=(y,z,c\cdot c'),
\]
where $c\cdot c'$ is the class of paths defined by glueing $c$ and $c'$
together.

\subsection{Constant Poisson structures}
This example combines the two previous ones.
Let us assume that $M=\bbR^n$ with a constant Poisson structure $\alpha$.
It is then possible to assume, if necessary after a linear change of 
coordinates, that $\alpha$ has the following block form:
\begin{gather*}
\alpha^{I\mu}=\alpha^{\mu\nu}=0,\qquad
I=1,\dots,r,\quad\mu,\nu=r+1,\dots,n,\\
(\alpha^{IJ})\text{ invertible},\qquad I,J=1,\dots,r,
\end{gather*}
where $r$ is the rank.

The the ``Gauss law'' reads
\[
\dd X^I+\alpha^{IJ}\,\eta_J=0,\qquad\dd X^\mu=0,
\]
and the infinitesimal symmetries are
\begin{align*}
\delta X^I &=\alpha^{IJ}\,\beta_J, & \delta X^\mu &=0,\\
\delta\eta_J&=\dd\beta_J, & \delta\eta_\mu&=\dd\beta_\mu.
\end{align*}

Thus, we can split $\bbR^n$ into $\bbR^r$ with symplectic structure
$(\alpha^{IJ})^{-1}$ and $\bbR^{n-r}$ with trivial Poisson structure.
By the previous two examples we get then
\[
\calG=\bbR^r\times\Bar\bbR^r\times\T^*\bbR^{n-r}
\]
with product
\[
(x,y,\xi,\pi)\sbullet(y,z,\xi,\pi')=(x,z,\xi,\pi+\pi').
\]

\subsection{The dual of a Lie algebra}\label{S-dLa}
Let $\fg^*$ be the dual of a Lie algebra $\fg$ with structure constants
in a given basis denoted by $f^{ij}_k$. 
The Kirillov--Kostant Poisson structure on $\fg^*$
is then given by the bivector field
\[
\alpha^{ij}(x)=f^{ij}_k\,x^k.
\]
In this case the Gauss law reads
\[
\dd X^i + f^{ij}_k\,X^k\,\eta_j=0,
\]
where $X$ is a map $I\to\fg^*$ and $\eta\in\Omega^1(I,\fg)$.

Let then $\beta$ be a map $I\to\fg$ that vanishes on the boundary of $I$.
The infinitesimal symmetries read
\[
\begin{aligned}
\delta X^i &=-f^{ij}_k\,X^k\,\beta_j,\\
\delta\eta_i &= \dd\beta_i+f_i^{jk}\,\eta_j\,\beta_k.
\end{aligned}
\]

We can rewrite the above equations in a more recognizable form if we
consider $\eta$ as the connection one-form for a $G$-bundle over $I$,
where $G$ is a Lie group whose Lie algebra is $\fg$.
The Gauss law becomes
\begin{equation}
\dd_\eta X = 0,\lab{detaX}
\end{equation}
while the infinitesimal symmetries now read
\begin{equation}
\begin{align}
\delta X &=\ad^*_{\beta} X,\\
\delta\eta &=\dd_\eta\beta,
\end{align}\lab{gaugesymm}
\end{equation}
that is, $\beta$ is an infinitesimal gauge transformation.

We define $\calG$ as the space of solutions of \eqref{detaX} modulo
gauge transformations connected to the identity.

We have then the following
\begin{Thm}
The phase space $\calG$ is diffeomorphic to $\T^*G$, where $G$
is the connected, simply connected Lie group whose Lie algebra is $\fg$.
The symplectic groupoid structure on $\T^*G$ is the one described in
Example~\ref{Liesg}.
\lab{thm-Glieg}
\end{Thm}

\subsubsection{Proof of Theorem~\ref{thm-Glieg}}
We first recall that $\T^*G$ is isomorphic to $\fg^*\times G$.
We then define a map 
\[
i\colon
\begin{array}[t]{ccc}
\calG &\to&\fg^*\times G\\
(X,\eta)&\mapsto&(X(0),\hol(\eta))
\end{array}
\]
where $0$ denotes the lower boundary of
$I$, and $\hol(\eta)$ is the parallel transport from the lower to the upper
boundary of $I$.

Next we want to define an inverse to $i$.
Let us then consider $(\xi,g)\in \fg^*\times G$. Since $G$ is connected,
there is a path $h:I\to G$ from the identity to $g$.
For such a path, we define
\[
\eta_{[h]} := h\dd h^{-1}.
\]
We then define
$X_{\xi,[h]}$ as the solution of \eqref{detaX} 
with initial condition $X_{\xi,[h]}(0)=\xi$ determined by $\eta_{[h]}$.
More precisely,
\[
X_{\xi,[h]} = \Ad^*_{h^{-1}} \xi.
\]
So $(X_{\xi,[h]},\eta_{[h]})$ is an element of $\calG$.
\begin{Lem}
Let $h$ and $l$ be two paths connecting the identity to the same 
element $g\in G$. Then $(X_{\xi,[h]},\eta_{[h]})=(X_{\xi,[l]},\eta_{[l]})$
in $\calG$.
\end{Lem}
\begin{Proof}
Let us consider the map
$\gamma:=hl^{-1}:I\to G$. 
Since $\gamma$ is the identity at the boundaries of $I$ and it
is in the connected component of the identity (as a consequence
of the fact that $G$ is simply connected), this is an allowed gauge
transformation.
Moreover, an easy computation proves that
\[
\begin{aligned}
(\eta_{[h]})^\gamma &= \gamma^{-1}\,\eta_{[h]}\,\gamma+
\gamma^{-1}\dd\gamma = \eta_{[l]},\\
(X_{\xi,[h]})^\gamma &= \Ad_\gamma^*X_{\xi,[h]}=X_{\xi,[l]}.
\end{aligned}
\]
Therefore, $(X_{\xi,[h]},\eta_{[h]})$ and $(X_{\xi,[l]},\eta_{[l]})$ define
the same element in $\calG$.
\end{Proof}
As a consequence we have the following well-defined map:
\[
j\colon
\begin{array}[t]{ccc}
\fg^*\times G&\to&\calG\\
(\xi,g)&\mapsto& (X_{\xi,g},\eta_g)
\end{array}
\]
with $(X_{\xi,g},\eta_g):=(X_{\xi,[h]},\eta_{[h]})$ for any 
path $h$ from the identity to $g$.

We then have the following
\begin{Lem}
The maps $i$ and $j$ are inverse to each other.
\end{Lem}
\begin{Proof}
Since $\eta_{[h]}$ is obtained from the trivial connection by
the gauge transformation $h^{-1}$ 
(which is not one of the symmetries we allow
since $h$ at the boundary is not the identity), we see immediately
that $\hol(\eta_{[h]})=g$. Moreover, $X_{\xi,g}(0)=\xi$ by definition.
So $i\circ j=\mathrm{id}$. 

Next we take a solution $(X,\eta)$ of \eqref{detaX}.
Let $l(u):=\hol^u(\eta)$ be the parallel transport from the lower 
boundary of $I$ till the point $u$. Notice that $l$ is a path
from the identity to $\hol(\eta)$.
Since moreover $\eta$ is equal to $l\dd l^{-1}$, we see that
$(X,\eta)=(X_{X(0),[l]},\eta_{[l]})$. But, from the previous
Lemma, we get then $(X,\eta)=(X_{X(0),\hol(\eta)},\eta_{\hol(\eta)})=
j\circ i(X,\eta)$.
\end{Proof}

To conclude the proof of Theorem~\ref{thm-Glieg}, we briefly discuss
the induced groupoid structure on $\fg^*\times G$.
Recalling that for us the the left and right projections correspond
to the boundary values $X(0)$ and $X(1)$, we obtain
\[
l(\xi,g)=\xi,\qquad r(\xi,g)=\Ad_{g^{-1}}^*\xi.
\]
The product is given by composition of solutions as in \eqref{compXeta}, 
and under this operation
the parallel transports also compose. So we get
\[
(\xi,g)\sbullet(\Ad_{g^{-1}}^*\xi,h)=(\xi,gh).
\]
After identifying $\fg^*\times G$ with $\T^*G$ by the map
described in Example~\ref{Liesg}, we recover the groupoid structure
described there.

\section{A singular phase space}\lab{S-6}

We want to discuss here an example proposed by Weinstein of a regular 
Poisson manifold that does not admit a symplectic groupoid and show what 
singularities arise in the construction of the phase space of the 
corresponding Poisson sigma model.

Let $M=\bbR^3\setminus\{0\}$ with Poisson bivector field
\[
\alpha^{ij}(x)=f(|x|)\,\pepsilon ijk\,x^k,\qquad f(R)\not=0\,\forall R>0.
\]
where $|\ |$ is the standard Euclidean norm.

For $f$ constant this Poisson manifold is equivalent to 
$\mathfrak{su}(2)\setminus \{0\}$ with
the Kirillov--Kostant Poisson structure, and the corresponding phase space
is $(\mathfrak{su}(2)^*\setminus \{0\})\times SU(2)$,
as described in the previous section.
If we introduce a non constant $f$, however, some problems may arise.

Observe first that, in any case, the symplectic leaves are the same as
in the case of $\mathfrak{su}(2)\setminus\{0\}$, i.e., spheres centered
at the origin. The symplectic form on these spheres is however rescaled by
a factor $f$, and the symplectic area $A$ of the sphere with radius $R$
is
\[
A(R) = \frac{4\,\pi\,R}{f(R)}.
\]
Then the observation of Weinstein \cite{W}, based on theorem of 
Dazord \cite{D}, is that
such a Poisson
manifold cannot have a symplectic groupoid if $A(R)$ is non constant
has critical points.

We want to see now how this condition arises in the construction of the phase
space. 

Namely, we have the following
\begin{Thm}
The phase space $\calG$ corresponding to $(M,\alpha)$ as above
is singular if{f} $A$ is non constant and has critical points.
\lab{thm-A}
\end{Thm}

\subsection{Proof}
In order to discuss the phase space $\calG$, 
it is convenient to use a vector notation; viz.,
we identify $(\bbR^3)^*$ and $\bbR^3$ using the Euclidean scalar product.
Moreover, we fix the volume form $du$ on the interval $I=[0,1]$.
Then both our fields $X$ and $\eta$ can be identified with
functions from $I$ to $\bbR^3$ that we denote by $\bX$ and $\bfeta$.
With these notations the Gauss law reads
\begin{equation}
\bX' + f(|\bX|)\,\bfeta\times\bX=0,
\lab{bX}
\end{equation}
where $\times$ denotes the cross product.

The infinitesimal symmetries can also be written in vector notation
after identifying $\beta$ with a map $\bbeta:I\to\bbR^3$:
\begin{equation}
\begin{split}
\dot\bX&=-f(|\bX|)\,\bbeta\times\bX,\\
\dot\bfeta&=\bbeta'+f(|\bX|)\,\bfeta\times\bbeta
+\frac{f'(|\bX|)}{|\bX|}\,(\bX\cdot\bfeta\times\bbeta)\,\bX,
\end{split}
\lab{dbX}
\end{equation}
where $\cdot$ is the Euclidean scalar product.

Given a map $\bv:I\to\bbR^3$ (e.g., $\bfeta$ or $\bbeta$), we define its
radial component $v_r$ and its tangential part $\bv_t$ w.r.t.\ $\bX$ by:
\begin{equation}
v_r(u):=\frac{\bv(u)\cdot\bX(u)}{|\bX(u)|},\qquad
\bv_t(u):=\bv(u)-v_r(u)\,\frac{\bX(u)}{|\bX(u)|}.
\lab{rt}
\end{equation}

Then we have the following:

\begin{Lem}
With the decomposition in \eqref{rt}, the Gauss law reads
\[
\bX' + f(|\bX|)\,\bfeta_t\times\bX=0,
\]
while the symmetries can be written as
\[
\begin{split}
\dot\bX&=-f(|\bX|)\,\bbeta_t\times\bX,\\
\dot\eta_r&=\beta_r'-\frac{f(|\bX|)}{|\bX|}\,(1-C(|\bX|))\,
(\bX\cdot\bfeta_t\times\bbeta_t),\\
\dot\bfeta_t&=\bbeta_t'+f(|\bX|)\,\bfeta_t\times\bbeta_t
+\frac{f(|\bX|)}{|\bX|^2}\,(\bX\cdot\bfeta_t\times\bbeta_t)\,\bX,
\end{split}
\]
with
\[
C(R) = \frac{R\,f'(R)}{f(R)} = 1-\frac{f(R)\,A'(R)}{4\,\pi}.
\]\lab{lem-rt}
\end{Lem}

\begin{Proof}
The Gauss law and the symmetry transformation for $\bX$ simply follow
from the fact that in a cross product or in a triple product
containing $\bX$ only tangential components of other vectors contribute.

For the symmetry transformation of $\bfeta$,
first of all we observe that $|\bX|'=|\bX|^{\cdot}=0$. Then we obtain
by \eqref{rt}, \eqref{bX} and \eqref{dbX}
the following identities:
\begin{align*}
\dot\eta_r&=(\dot\bfeta)_r-\frac{f(|\bX|)}{|\bX|}\,
\bX\cdot\bfeta_t\times\bbeta_t,\\
\beta_r'&=(\bbeta')_r+\frac{f(|\bX|)}{|\bX|}\,
\bX\cdot\bfeta_t\times\bbeta_t.
\end{align*}
These yield immediately the symmetry equation for $\eta_r$.

To obtain the symmetry equation for $\bfeta_t$, we first observe that
\[
\bbeta_t'=\bbeta'-\beta_r'\frac\bX{|\bX|}+\frac{f(|\bX|)}{|\bX|}\,\beta_r\,
\bfeta_t\times\bX.
\]
Then we get
\begin{multline*}
\dot\bfeta_t=\dot\bfeta-\dot\eta_r\frac\bX{|\bX|}
+\frac{f(|\bX|)}{|\bX|}\,\eta_r\,
\bbeta_t\times\bX=\\
=\bbeta_t'+f(|\bX|)\,\bfeta\times\bbeta-\frac{f(|\bX|)}{|\bX|}\,\beta_r\,
\bfeta_t\times\bX+\frac{f(|\bX|)}{|\bX|}\,\eta_r\,
\bbeta_t\times\bX+\\ 
+\frac{f(|\bX|)}{|\bX|^2}\,(\bX\cdot\bfeta_t\times\bbeta_t)\,\bX,
\end{multline*}
which, after using again \eqref{rt}, leads to the desired identity.
\end{Proof}

Observe now that the original case of $\mathfrak{su}(2)$ 
is recovered by setting
$f\equiv1$ and $C\equiv0$ in the equations displayed in Lemma~\ref{lem-rt}.
On the other hand, the critical case $A'(R)=0$ corresponds to $C(R)=1$.

Let us begin considering solutions with $C(|\bX|)\not=1$. In this case,
we can define new variables as follows:
\begin{align}
a_r &= \frac{f(|\bX|)}{1-C(|\bX|)}\,\eta_r, &
\ba_t &= f(|\bX|)\,\bfeta_t,\lab{resca}\\
b_r &= \frac{f(|\bX|)}{1-C(|\bX|)}\,\beta_r, &
\bb_t &= f(|\bX|)\,\bbeta_t.
\end{align}
Then we obtain the Gauss law in the form
\[
\bX' + \ba_t\times\bX=0,
\]
while the symmetries read now
\[
\begin{split}
\dot\bX&=-\bb_t\times\bX,\\
\dot a_r&=b_r'-\frac{1}{|\bX|}\,
(\bX\cdot\ba_t\times\bb_t),\\
\dot\ba_t&=\bb_t'+\ba_t\times\bb_t
+\frac{1}{|\bX|^2}\,(\bX\cdot\ba_t\times\bb_t)\,\bX.
\end{split}
\]
Thus we have recovered, in the new variables, the case of $\mathfrak{su}(2)$.
Proceeding now as in the proof of Theorem~\ref{thm-Glieg} (namely, taking
holonomies of $\ba$ as coordinates), we find that 
the fiber of the left projection of $\calG$                      
 over a point $\bx\in M$
with $C(|\bx|)\not=1$ is diffeomorphic to $SU(2)$.

On the other hand, when $C(|\bx|)=1$, the above change of variables is not 
defined.
In this case we may however choose the following complete set of
invariant functions:
\[
\bx:=\bX(0),\qquad
\by:=\frac{\bX(1)}{|\bx|},\qquad 
\pi:=\int_I\eta_r.
\]
That is, the fiber of the left projection
of $\calG$ over $\bx$ 
with $C(|\bx|)=1$ is diffeomorphic to $S^2\times\bbR$.

If $C\equiv1$---i.e., if $A$ is constant---then $\calG$ is the
smooth manifold $\bbR^+\times S^2\times S^2\times\bbR$.

To better visualize the singularities in the general case,
let us pick up a neighborhood $U$ of a point in $\bbR^+$ where
$A'$ vanishes but $A$ is non constant. Let $V$ be a neighborhood of
a point in $S^2$. We want to show that 
$\calG_{UV}:=l^{-1}(U\times V)$ is not a manifold.
We can describe $\calG_{UV}$ as follows. Given a solution $(\bX,\bfeta)$,
we can always transform it into a solution with $\eta_r$ constant (just    
take a transformation generated by $\bbeta$ with $\bbeta_t=0$). Under small
gauge transformations such a solution is characterized by the values of $\bX$
at the endpoints and the value of $\eta_r$. If $C(|\bX|)=1$, there is no
way of changing $\eta_r$ into another constant. On the other hand,
if $C(|\bX|)\not=1$, large gauge transformations can send $\eta_r$ into
another constant that differs from the former by a multiple of
$4\pi\,[1-C(R)]/f(R)$ (we are trivializing the Hopf bundle $SU(2)\to S^2$
over $V$ taking into account the rescaling \eqref{resca}). 
Therefore, $\calG_{UV}=V\times V\times Q$, where $Q$ is
the quotient of
$U\times\bbR$ by the equivalence relation
\[
(R,p)\sim \left(R,p+\frac{4\pi\,[1-C(R)]}{f(R)}\right).
\]

\section{The phase space of the Poisson sigma model with 
two-dimensional target}\lab{S-7}
Let $U$ be a domain in $\bbR^2$ with a given Poisson bivector field 
$\alpha$. After choosing coordinates, it is always possible to write
\[
\alpha^{ij}(x^1,x^2) = \epsilon^{ij}\,\phi(x^1,x^2),
\quad\phi\in C^\infty(U).
\]
We also fix the volume form $du$ on $I$ and then identify $\eta$
with a map $I\to\bbR^2$.

With these choices, the ``Gauss law''simply reads
\begin{equation}
\begin{aligned}
(X^1)' + \phi(X^1,X^2)\,\eta_2 &=0,\\
(X^2)' - \phi(X^1,X^2)\,\eta_1 &=0.
\end{aligned}
\lab{Gl}
\end{equation}

The infinitesimal symmetries read then
\begin{equation}
\begin{aligned}
\delta X^1&=-\phi\,\beta_2,\\
\delta X^2&=\phi\,\beta_1,\\ 
\delta\eta_1&=\beta_1'+\de_1\phi\,(\eta_1\,\beta_2-\eta_2\,\beta_1),\\
\delta\eta_2&=\beta_2'+\de_2\phi\,(\eta_1\,\beta_2-\eta_2\,\beta_1),
\end{aligned}
\lab{bsymm}
\end{equation}
where $\de_i\phi$ is a shorthand notation for $\de\phi/\de x^i$, and the 
infinitesimal generators $\beta_i$ are continuously
differentiable maps $I\to \bbR^2$
with the conditions
\begin{equation}
\beta_i(0)=\beta_i(1)=0,\quad i=1,2.
\lab{b0}
\end{equation}

We will denote by  $\tcalG$  the phase space of solutions of \eqref{Gl} 
modulo 
the symmetries generated by \eqref{bsymm}. More precisely, we
first introduce the Banach spaces
$C^0(I,\bbR^2)$, $C^1(I,\bbR^2)$ and $C^1_0(I,\bbR^2)$. Then we consider
the Banach manifold $C^1(I,U)$ modeled on $C^1(I,\bbR^2)$.
With these notations we can finally write
\[
\tcalG :=\frac{\{
(X,\eta)\in C^1(I,U)\times C^0(I,\bbR^2)\ |\ (X,\eta) \text{ satisfies \eqref{Gl}}
\}}
{\{
\text{symmetries \eqref{bsymm} with $\beta\in C^1_0(I,\bbR^2)$}
\}}.
\]

In the rest of this section we will study $\tcalG$.
Namely, in subsection~\ref{ssec-ttG} we will give an equivalent
but easier description
of $\tcalG$, and in subsection~\ref{ssec-G} we will
show that the latter is diffeomorphic to a submanifold $\calG$ of $\bbR^4$, at
least if all the symplectic leaves of $U$ are simply connected; 
in subsection~\ref{ssec-Gg}
we will describe the product structure for $\calG$ induced from the
composition of paths $X:I\to U$; finally, in subsection~\ref{ssec-Gs} 
we will derive the symplectic structure for $\calG$
from the symplectic structure on the space of fields
$(X,\eta)$.

\subsection{An equivalent description of the phase space}\lab{ssec-ttG}
>From now on, by abuse of notation, we will write $\phi$ for $\phi\circ X$.

The Gauss law \eqref{Gl} implies an equation for $\phi$, viz.,
\begin{equation}
\phi'= T\,\phi,
\lab{phidot}
\end{equation}
with
\begin{equation}
T := \de_2\phi\,\eta_1 - \de_1\phi\,\eta_2.
\lab{defT}
\end{equation}
The solution of \eqref{phidot} is simply given by
\begin{equation}
\phi(X^1(u),X^2(u))=\phi_0\,H(u),
\lab{phisol}
\end{equation}
where $\phi_0$ is a shorthand notation for $\phi(X^1(0),X^2(0))$ and
\begin{equation}
H(u):=\exp\int_0^u T(v)\,dv.
\lab{defH}
\end{equation}

It is then useful to define the following change of variables:
\begin{equation}
E_i := \eta_i\,H,\lab{defE}
\end{equation}
Notice that the map $(X,\eta)\mapsto(X,E)$
is a smooth map from $C^1(I,U)\times C^0(I,\bbR^2)$ into itself.

With these new variables, we can rewrite the Gauss law \eqref{Gl} as
\begin{equation}
\begin{aligned}
(X^1)' + \phi_0\,E_2 &=0,\\
(X^2)' - \phi_0\,E_1 &=0.
\lab{mGl}
\end{aligned}
\end{equation}
Notice that every solution of \eqref{Gl} determines a solution of
\eqref{mGl} via \eqref{defE}. The converse, however, is not true in general.

Assume in fact that $(X,E)$ is a solution of \eqref{mGl}.
Then we get the following equation for $\phi$:
\begin{equation}
\phi'= \tT\,\phi_0,
\lab{mphidot}
\end{equation}
with
\begin{equation}
\tT := \de_2\phi\,E_1 - \de_1\phi\,E_2.
\lab{deftT}
\end{equation}
Comparing the solution of \eqref{mphidot} with \eqref{phisol}, we get
\begin{equation}
H(u) = 1 + \int_0^u \tT(v)\,dv.
\lab{Ht}
\end{equation}
By comparison with \eqref{defH}, we see that a solution $(X,E)$ of 
\eqref{mGl} determines a solution $(X,\eta)$ of \eqref{Gl} if{f}
the following condition is satisfied:
\begin{equation}
H(u)>0,\quad\forall u\in I.
\lab{cond}
\end{equation}
So we have the following
\begin{Lem}
Solutions of \eqref{Gl} are mapped by \eqref{defE} into solutions of
\eqref{mGl} satisfying \eqref{cond} and vice versa.
\end{Lem}
As for the symmetries acting on $(X,E)$, we introduce
\begin{equation}
e_i := \beta_i\,H.
\lab{defe}
\end{equation}
Observe here that the map $(X,\eta,\beta)\mapsto e$ is a smooth
map from $C^1(I,U)\times C^0(I,\bbR^2)\times C^1_0(I,\bbR^2)\to C^1_0(I,\bbR^2)$.

Then we have the following:
\begin{Lem}
Under the infinitesimal symmetry \eqref{bsymm}, the variables $(X,E)$
defined via \eqref{defE}
in terms of a solution $(X,\eta)$ of \eqref{Gl}
change as follows:
\begin{equation}
\begin{aligned}
\delta X^1&=-\phi_0\,e_2,\\
\delta X^2&=\phi_0\,e_1,\\ 
\delta E_1&=e_1',\\
\delta E_2&=e_2'.
\end{aligned}
\lab{mbsymm}
\end{equation}
Conversely, if $(X,E)$ is a solution of \eqref{mGl} satisfying to
\eqref{cond}, then the infinitesimal symmetry \eqref{mbsymm}
implies the infinitesimal symmetry \eqref{bsymm} on 
the variables $(X,\eta)$ obtained by inverting \eqref{defE}.
\end{Lem}
\begin{Proof}
The first two equations are immediately obtained by \eqref{phisol}
and \eqref{defe}.

As for the two other equations, we first observe that
\[
\delta T = \tau',
\]
with
\[
\tau := \de_2\phi\,\beta_1 - \de_1\phi\,\beta_2.
\]
In fact,
\begin{multline*}
\delta T =
\delta(\epsilon^{ij}\, \eta_i\,\de_j\phi)=\\
=\epsilon^{ij}\,\beta_i'\,\de_j\phi 
+\epsilon^{ij}\,\de_i\phi\epsilon^{kl}\,\eta_k\,\beta_l\,\de_j\phi
-\epsilon^{ij}\,\eta_i\,\phi_{jk}\,\epsilon^{kl}\,\phi\,\beta_l=\\
=\frac{\dd}{\dd u}(\epsilon^{ij}\, \beta_i\,\de_j\phi)=
\tau',
\end{multline*}
where we have made use of \eqref{Gl} and \eqref{bsymm}.
{}From this we get
\[
\delta H = \tau H.
\]
Finally,
\[
\delta E_i =
\delta\eta_i\,H + \eta_i\,\tau\,H =
\frac{\dd}{\dd u}(\beta_i\,H) -
(\beta_i\,T-\eta_i\,\tau-\de_i\phi\,\epsilon^{kl}\,\eta_k\,\beta_l)\,H.
\]
A direct computation shows that the 
terms in the second brackets cancel, 
so we have proved the first part of the Lemma.

As for the second part, we observe that
\[
\delta\tT = {\Hat{\tau}}',
\]
with $\Hat\tau=\de_2\phi\,e_1-\de_1\phi\,e_2$. As a consequence,
$\delta H = \Hat\tau$. Observing then that $\tT=H\,T$ and 
$\Hat\tau=H\,\tau$, we get
\begin{multline*}
\delta\eta_i = \delta\pfrac{E_i}H=
\frac{\delta E_i}H-\frac{E_i}{H^2}\,\delta H=
\frac{e_i'}{H}-\frac{E_i}{H^2}\,\Hat\tau=\\
=\frac{\frac{\dd}{\dd u}(\beta_i\,H)}H-\eta_i\,\tau=
\beta_i'+\beta_i\,T-\eta_i\,\tau,
\end{multline*}
from which \eqref{bsymm} follows.
\end{Proof}

We now define a new phase space:
\[
\ttcalG :=\frac{\{
(X,E)\in C^1(I,U)\times C^0(I,\bbR^2)\ |\ (X,E) \text{ satisfies \eqref{mGl}
and \eqref{cond}}
\}}
{\{
\text{symmetries \eqref{mbsymm} with $e\in C^1_0(I,\bbR^2)$}
\}}.
\]
Then the preceding discussion, and in particular the two Lemmata,
prove the following
\begin{Prop}
If $\ttcalG$ is a smooth manifold, then $\tcalG$ and $\ttcalG$ are
diffeomorphic.\lab{prop-ttcalG}
\end{Prop}
In the next subsection we will prove that $\ttcalG$ is actually a
smooth $4$-manifold, at least under the following
\begin{Ass}
We assume that all the symplectic leaves of $(U,\alpha)$ are simply
connected.\lab{ass-sim}
\end{Ass}
Observe that for example $\bbR^2$ with $\phi=(x^1)^2+(x^2)^2$ will not
be allowed.

\subsection{The phase space is a smooth manifold}\lab{ssec-G}
We begin by defining some invariants of $\ttcalG$.
The first are the initial conditions of $X$, viz.,
\begin{equation}
x^i := X^i(0);
\lab{x0}
\end{equation}
the others are the following integrals:
\begin{equation}
\pi_i := \int_0^1 E_i(u)\,du.
\lab{pi}
\end{equation}
All of them are invariant under \eqref{mbsymm} since $e\in C^1_0(I,\bbR^2)$.

In this way we get a well-defined, smooth map $i:\ttcalG\to U\times\bbR^2$.
This map is however not surjective because of \eqref{cond}. 

We want then to define an appropriate domain in $\bbR^4$
so that $i$ becomes a diffeomorphism. 

We first define the continuous map $x_f:U\times\bbR^2\to\bbR^2$ by
\begin{equation}
x^i_f = x^i-\phi(x^1,x^2)\,\epsilon^{ij}\,\pi_j,
\lab{xf}
\end{equation}
and then
\[
V:=\{p\in U\times\bbR^2\ \vert\  x_f(p)\in U\}.
\]
Remark that $x_f$ can also be seen
as the final point of a solution $X$ of \eqref{mGl}, with $x$ and $\pi$
given by \eqref{x0} and \eqref{pi}.

Next we define the map $h:V\to\bbR$ by
\begin{equation}
h(x^1,x^2,\pi_1,\pi_2) := 
\begin{cases}
\dfrac{\phi(x^1_f,x^2_f)}
{\phi(x^1,x^2)} & \text{if $\phi(x^1,x^2)\not=0$},\\
1-\pi_2\,\de_1\phi(x^1,x^2)+\pi_1\,\de_2\phi(x^1,x^2)
 & \text{if $\phi(x^1,x^2)=0$}.
\end{cases}\lab{defh}
\end{equation}
Then we define 
\[
\boldsymbol{\ucalG}:=
\{p\in V\ \vert\ h(p)>0\},
\]
and finally we denote by
${\calG}$ the connected component of $\bucalG$
containing $U\times\{(0,0)\}$.
\begin{Lem}
${\calG}$ is a 4-manifold.
\lab{uG4mfld}
\end{Lem}
\begin{Proof}
We just have to prove that $h$ is continuous.
To do this, we observe that the two definitions for $h$ are continuous
when restricted to the appropriate subset.

Since the zero locus of $\phi$ is closed, we only have to check
that, for any sequence in the complement of the zero locus that approaches
a point in the zero locus, the limit of the first expression yields
the second expression. This can be easily proved by Taylor expanding
the numerator.

To prove that the connected component we are interested in is not empty,
it is enough to observe that $h(x^1,x^2,0,0)=1$, $\forall(x^1,x^2)\in U$.
\end{Proof}

\begin{Exa}[Semiclassical quantum plane]\lab{exa-qp}
Let $U=\bbR^2$ and $\phi(x^1,x^2)=x^1\,x^2$.
In this case the map $h$ simply reads
\[
h(x^1,x^2,\pi_1,\pi_2) = (1-x^2\pi_2)(1+x^1\pi_1).
\]
In the connected component of $h^{-1}(\bbR^+)$ both factors
must be positive. So we get
\[
\calG=
\{
(x^1,x^2,\pi_1,\pi_2)\in\bbR^4\ \vert\ x^1\pi_1>-1,\ x^2\pi_2<1
\}.
\]
The 2-dimensional symplectic leaves are the four open quadrants.
Over each point $(x^1,x^2)$
of one of these leaves, the fiber is given by those vectors
$(\pi_1,\pi_2)$ such that the linear trajectory with constant velocity
\[
(-\phi(x^1,x^2)\,\pi_2,\phi(x^1,x^2)\,\pi_1)
\]
is entirely contained
in the same symplectic leaf for all times $t\le1$. 
Over points in the zero locus of $\phi$, i.e., the axes, the fiber is
the whole of $\bbR^2$, for the velocity here is zero.

Observe that this simple description of $\calG$ is possible 
whenever all symplectic
leaves in $U$ are convex; e.g., when $U=\bbR^2$ and 
$\phi(x^1,x^2)=(x^1)^r\,(x^2)^s$, $r,s\ge0$.
\end{Exa}


The central result of this section is the following:
\begin{Thm}
Under Assumption~\ref{ass-sim},
the phase spaces $\tcalG$ and $\ttcalG$ are diffeomorphic
to $\calG$.
\lab{GGG}
\end{Thm}






\subsection{Proof of Theorem~\ref{GGG}}
In view of Proposition~\ref{prop-ttcalG}, we have
only to prove that
$\calG$ is diffeomorphic to $\ttcalG$.

The idea is to show that the mapping given by \eqref{x0} and 
\eqref{pi} defines a diffeomorphism.

\begin{Lem}
There is a well-defined smooth map
\[
i:\begin{array}[t]{ccc}
\ttcalG &\to& \calG \\{}
[(X,E)] &\mapsto& \left({X(0),\int_0^1 E(u)\,du}\right)
\end{array}
\]
\end{Lem}
\begin{Proof}
Since \eqref{x0} and \eqref{pi} are invariant under the symmetries
\eqref{mbsymm}, this map descends to $\ttcalG$.

We want then to show that the image of this map is contained in $\calG$.

First of all, we observe that $x_f=X(1)$; so automatically
\[
x_f(x^1,x^2,\pi_1,\pi_2)\in U.
\]

Then we want to show that \eqref{cond} implies $h>0$.
Consider first the case when $\phi_0\not=0$. Then, by \eqref{phisol},
$\phi(u)/\phi_0=H(u)>0,\forall u\in I$. In particular, for $u=1$, this
implies $\phi(x_f)/\phi(x)>0$.
The other possibility is when $\phi_0=0$. In this case $X$ is constant,
and by \eqref{Ht} we get that $H(1)=h$.

Since for every solution $(X,E)$ of \eqref{mGl} and \eqref{cond}
there is a solution
with the same initial condition and $\Tilde{E}=E/\lambda$, $\lambda>0$, then
$\ttcalG$ is connected, so its image is contained in a connected
component of the set $h^{-1}(\bbR^+)$. The existence of constant
solutions with $E=0$ implies that this connected component contains
$U\times\{(0,0)\}$.

So we have proved that $[(X,E)]\mapsto(x,\pi)$ is a well-defined
map from $\ttcalG$ to $\ucalG$.
\end{Proof}

We now want to define an inverse map 
\[
j:\begin{array}[t]{ccc}
\calG&\to&\ttcalG\\
(x,\pi)&\mapsto&[(X,E)]
\end{array}
\]
We consider two cases:
\begin{enumerate}
\item $\phi(x)\not=0$. 
We take $X$ equal to
any path joining $x$ to $x_f(x,\pi)$
that is completely contained in the symplectic
leaf.
Then we set $E_i(u)=\epsilon_{ij}\,(X^j)'(u)/\phi(x)$.
\item $\phi(x)=0$.
We set $X(u)=x$ and $E(u)=\pi$ $\forall u\in I$.
\end{enumerate}

It is not difficult to see that the image of $j$ is a solution
of \eqref{mGl} and \eqref{cond}.
\begin{Lem}
Let $(X,E)$ and $(\Tilde{X},\Tilde{E})$ be two solutions determined as above.
Then, under Assumption~\ref{ass-sim}, 
they define the same element in $\ttcalG$.
\end{Lem}
\begin{Proof}
In the case when $\phi(x)=0$, we completely specified the
solution; so $(X,E)=(\Tilde{X},\Tilde{E})$.

Let us consider then the case $\phi(x)\not=0$.
Since any symplectic leaf is simply connected by Assumption~\ref{ass-sim},
there is a path $X(u,s)$ connecting $X$ to $\Tilde{X}$.
More precisely, $X(u,0)=X(u)$, $X(u,1)=\Tilde{X}(u)$,
$X(0,s)=x$, $X(1,s)=x_f$, and $X(\sbullet,s)$ is entirely
contained in the symplectic leaf of $x$.
We set then
\[
e_i=\epsilon_{ij}\,\frac{\dot X^j}{\phi(x)}
\]
and integrate the infinitesimal symmetry \eqref{mbsymm} obtaining
\[
E_i(u,s)=E_i(u)+\int_0^s e_i'(u,\sigma)\,d\sigma.
\]
\end{Proof}

As a consequence the map $j$ is well-defined, and it is not difficult
to prove that it is smooth. Moreover, we have the following
\begin{Lem}
Under Assumption~\ref{ass-sim},
the maps $i$ and $j$ are inverse to each other.
\end{Lem}
\begin{Proof}
The identity $i\circ j=\mathrm{id}$ is trivial. 

We want to prove that also $j\circ i=\mathrm{id}$. 

Let us begin with the case $\phi_0\not=0$. In this case $j\circ i[(X,E)]$
is a solution $(\Tilde{X},\Tilde{E})$ so that $\Tilde{X}$ has the same 
end-points of $X$. Thus, as in the proof of the previous Lemma,
there is a symmetry that relates them.

In the case when $\phi_0=0$, we must prove that, given a solution
$(X,E)$, then $(\Tilde{X},\Tilde{E}):=j\circ i(X,E)$ is an equivalent solution.

First we observe that $X=\Tilde{X}$ since both are constant solutions with
the same starting point. 

We have then to find a symmetry that sends $E$ to $\Tilde{E}$.
To do so, we define the following element of $C^1_0(I,\bbR^2)$:
\[
e_i(u):=\int_0^u (\Tilde{E}_i-E_i(v))\,dv.
\]
Then we consider the path in $C^0(I,\bbR^2)$ given by
\[
E_i(u,s):=s\,\Tilde{E}_i+(1-s)\,E_i(u),\qquad s\in[0,1].
\]
We have then
\[
\begin{aligned}
E_i(u,0)&=E_i(u),\\
E_i(u,1)&=\Tilde{E}_i,\\
\dot E_i(u,s) &= e_i'(u).
\end{aligned}
\]
So we can go from $E$ to $\Tilde{E}$ via a symmetry transformation 
\eqref{mbsymm}.

Since $\phi_0=0$, the corresponding path of paths $X(u,s)$ is constant
and equal to $X(0)$.

To complete the proof,
we have only to check that condition \eqref{cond} is satisfied for
any intermediate value $u\in[0,1]$.
To do this we just observe that, by definition,
\begin{multline*}
H(u,s)=\\
=1+\de_2\phi\,[s\,\Tilde{E}_1+(1-s)E_1(u)]
-\de_1\phi\,[s\,\Tilde{E}_2+(1-s)E_2(u)]=\\
=s\,A+(1-s)B(u),
\end{multline*}
with $A=1+\de_2\phi\,\Tilde{E}_1-\de_1\phi\,\Tilde{E}_2$ and
$B(u)=1+\de_2\phi\,E_1(u)-\de_1\phi\,E_2(u)$.
Since $A>0$ and $B(u)>0\ \forall u\in I$, we get that
$H(u,s)>0\ \forall(u,s)\in I\times[0,1]$.
\end{Proof}
This concludes the proof of Theorem \eqref{GGG}.

\subsection{The product on $\calG$}\lab{ssec-Gg}
We will describe $\calG$ in terms of local coordinates
$(x,\pi)$ with $x\in U$, $\pi\in\bbR^2$.

We define the two projections $r,l:\calG\to U$, by
\[
l(x,\pi):=x,\qquad
r(x,\pi):=x_f(x,\pi),
\]
which correspond to the initial and final point of the given
solution in $\tcalG$ as prescribed by \eqref{lrXeta}.

Let us consider now another point $(\tilde x,\tilde\pi_1)\in\calG$,
with $\tilde x=x_f(x,\pi)$.

Then we look for the solutions $(X,\eta)$ and $(\Tilde X,\Tilde\eta)$
in $\tcalG$ that correspond to the points in $\calG$ 
described above.
In particular we choose the solutions so that the tangent at $X(1)$
is equal to the tangent at $\Tilde X(0)$.
So we can compose the solutions in a differentiable way
as in \eqref{compXeta}.
We now want to compute the point $(\Hat x,\Hat\pi)\in\calG$
corresponding to the new solution $(\Hat X,\Hat\eta)$.
We immediately get
\begin{equation}
\Hat x= x.
\lab{prodxc}
\end{equation}

As for $\Hat\pi$, we use \eqref{pi} and \eqref{defE} obtaining
\begin{multline*}
\Hat\pi=
\int_0^1\Hat H(u)\,\Hat\eta(u)\,du =
\int_0^1\Hat H\left(\frac u2\right)\,\eta(u)\,du +
\int_0^1\Hat H\left(\frac {u+1}2\right)\,\Tilde\eta(u)\,du.
\end{multline*}
By \eqref{defH} we get then
\[
\Hat H\left(\frac u2\right)=H(u),\qquad
 H\left(\frac {u+1}2\right)=H(1)\,\Tilde H(u),
\]
for $u\in[0,1]$. Thus, we get
\begin{equation}
\Hat\pi = \pi + h(x,\pi)\,\Tilde\pi.
\lab{prodpi}
\end{equation}

\subsection{The symplectic structure on $\calG$}\lab{ssec-Gs}
As in the general description, we consider the constant symplectic structure
on $C^1(I,U)\times C^0(I,\bbR^2)$ determined by the action $\int_0^1 \eta_i(u)\,(X^i)'(u)$;
viz.,
\[
\omega((\alpha,\beta),(\Tilde{\alpha},\Tilde{\beta})):=
\int_0^1 [\Tilde{\alpha}(u)\,\beta(u)-\alpha(u)\,\Tilde{\beta}(u)]\,du,
\]
for $(\alpha,\beta),(\Tilde{\alpha},\Tilde{\beta})\in \T(C^1(I,U)\times C^0(I,\bbR^2))$.

In order to perform the computations of this subsection, it is however more 
convenient to work with the corresponding Poisson structure that
we write
\[
\Poiss{\eta_i(u)}{X^j(v)}=\delta_i^j\,\delta(u-v),
\]
while all other brackets vanish. As usual in infinite dimensional cases,
the Poisson bracket is defined only for a certain class of functions.

We now want to determine the induced Poisson structure on $\calG$.

By the general argument we get
\begin{equation}
\Poiss{x^1}{x^2}=\phi(x^1,x^2).
\lab{Pxx}
\end{equation}
An easy computation yields
\[
\Poiss{T(u)}{X^i(0)}=\epsilon^{ij}\,\de_j\phi_0\,\delta(u).
\]
Thus we get
\begin{multline}
\Poiss{x^i}{\pi_j}=
\Poiss{X^i(0)}{\int_0^1 H(u)\,\eta_j(u)\,du}=\\
=-\delta_i^j+\int_0^1 H(u)\,\eta_j(u)\,\Poiss{X^i(0)}{\int_0^u T(v)\,dv}=\\
=-\delta_i^j-\pi_j\,\epsilon^{ik}\,\de_k\phi(x^1,x^2).
\lab{Pxpi}
\end{multline}
Finally, we have the most complicated bracket, that is, 
$\Poiss{\pi_1}{\pi_2}$.
\begin{Lem}
Let us consider the function $\psi:\ucalG\to\bbR$ defined by
\begin{equation}
\psi:=
\begin{cases}
\dfrac{1+\pi_1\,\de_2\phi-\pi_2\,\de_1\phi
-h}{\phi}
&\text{if $\phi\not=0$},\\
\pi_1\,\pi_2\,\de_1\de_2\phi-\frac12\,(\pi_1)^2\,(\de_2)^2\phi-
\frac12\,(\pi_2)^2\,(\de_1)^2\phi
&\text{if $\phi=0$},
\end{cases}\lab{defpsi}
\end{equation}
with $h$ defined in \eqref{defh}.
Then $\psi$ is smooth and
\begin{equation}
\Poiss{\pi_1}{\pi_2}=\psi(x^1,x^2,\pi_1,\pi_2).
\lab{Ppipi}
\end{equation}
\end{Lem}
\begin{Proof}
The smoothness of $\psi$ is proved by Taylor expanding $h$ in the first
expression.

As for the second assertion, we first observe the following 
useful identities:
\begin{align*}
\Poiss{\eta_i(u)}{T(v)} &=
\delta(u-v)\,\epsilon^{kl}\,\eta_k(u)\,\de_i\de_l\phi(u),\\
\Poiss{T(u)}{T(v)} &=0,\\
\intertext{which imply}
\Poiss{\eta_i(u)}{H(v)} &=
\theta(v-u)\,\epsilon^{kl}\,H(v)\,\eta_k(u)\,\de_i\de_l\phi(u),\\
\Poiss{H(u)}{H(v)} &=0.
\end{align*}
Then a straightforward computation yields
\begin{multline*}
\Poiss{\pi_1}{\pi_2}=
\int_0^1du\int_0^1dv\ [
E_1(u)\,E_2(v)\,(\theta(u-v)\,
\de_1\de_2\phi_v+\theta(v-u)\,\de_1\de_2\phi_u)+\\
-\theta(u-v)\,E_1(u)\,E_1(v)\,(\de_2)^2\phi_v
-\theta(v-u)\,E_2(u)\,E_2(v)\,(\de_1)^2\phi_v
],
\end{multline*}
where $\phi_v$ is a short-hand notation for $\phi(X^1(v),X^2(v))$.

In the case when $\phi_0=0$, the solution $X$ is constant.
So we can take all the terms of the form $\de_i\de_j\phi$ out of
the integral. What is left, thanks to \eqref{pi}, yields
the second formula in \eqref{defpsi}.

If however $\phi_0\not=0$, we multiply both sides by $(\phi_0)^2$
and then use \eqref{mGl}, obtaining
\begin{multline*}
(\phi_0)^2\,\Poiss{\pi_1}{\pi_2}=
\int_0^1du\int_0^1dv\ \Big[
\theta(u-v)\,(X^2)'(u)\,\frac\dd{\dd v}\de_2\phi_v+\\
+\theta(v-u)\,(X^1)'(v)\,\frac\dd{\dd u}\de_1\phi_u
\Big].
\end{multline*}
A simple integration yields then
\[
(\phi_0)^2\,\Poiss{\pi_1}{\pi_2}=
(X^2(1)-X^2(0))\,\de_2\phi_0+
(X^1(1)-X^1(0))\,\de_1\phi_0
-\phi_1+\phi_0,
\]
which is the first formula in \eqref{defpsi} thanks to \eqref{x0}, 
to \eqref{xf} and to the identity $h=\phi_1/\phi_0$.
\end{Proof}

The brackets of the coordinates define a bivector field on $\ucalG$, which
we will denote by $\frP$, through the relation
\begin{equation}
\Poiss fg = \frP(\dd f,\dd g).
\lab{PfrP}
\end{equation}
Locally, in the basis corresponding
to the coordinates $x^1,x^2,\pi_1,\pi_2$, we write this bivector
field in matrix form as
\begin{equation}
\matP=
\begin{pmatrix}
0 & \phi & -1-\pi_1\,\de_2\phi & -\pi_2\,\de_2\phi \\
-\phi & 0& \pi_1\,\de_1\phi & -1+\pi_2\,\de_1\phi \\
1+\pi_1\,\de_2\phi & -\pi_1\,\de_1\phi & 0 & \psi \\
\pi_2\,\de_2\phi & 1-\pi_2\,\de_1\phi & -\psi & 0
\end{pmatrix}.
\lab{defP}
\end{equation}
This matrix is always invertible thanks to the condition $h>0$.
We will exhibit the corresponding 2-form 
$\omega_\calG$ in the next subsection.

\subsection{Summary}
We started with a 2-dimensional domain $U$ with a bivector field
$\alpha^{ij}=\epsilon^{ij}\,\phi$, $\phi\in C^\infty(U)$, such that
Assumption~\ref{ass-sim} holds.

Using $\phi$ we defined the 4-dimensional domain $\ucalG$ as at the beginning
of subsection~\ref{ssec-G}, viz., as the connected component containing 
$U\times\{(0,0)\}$ of the set
\[
\{
(x,\pi)\in U\times\bbR^2\ |\ x_f(x,\pi)\in U,\ h(x,\pi)>0
\},
\]
with $x_f$ and $h$ defined in \eqref{xf} and \eqref{defh}.

Next we obtained the left and right projections $l,r:\ucalG\to U$ by
\begin{align*}
l(x,\pi) &= x,\\
r(x,\pi) &= x_f(x,\pi) = x-\alpha\,\pi.
\end{align*}

Given two points $(x,\pi)$ and $(\tilde x, \tilde\pi)$ in $\ucalG$ with
$r(x)=l(\tilde x)$, we got the product
\begin{equation*}
(x,\pi)\sbullet (\tilde x, \tilde\pi) =
(x, \pi + h(x,\pi)\,\tilde\pi).
\end{equation*}

Finally, in \eqref{defP} we defined a bivector field $\frP$
whose inverse exists and is given by the following 2-form:
\begin{multline*}
\omega_\calG=[
\psi\, \dd x^1\wedge\dd x^2 + (1-\pi_2\,\de_1\phi)\,\dd x^1\wedge\dd\pi_1 +
\pi_1\,\de_1\phi\,\dd x^2\wedge\dd\pi_2+\\
 -\pi_2\,\de_2\phi\,\dd x^2\wedge\dd\pi_1 
+(1+\pi_1\,\de_2\phi)\,\dd x^2\wedge\dd\pi_2 -\phi\,\dd\pi_1\wedge\dd\pi_2
]/h.
\end{multline*} 

{}From the general results of Section~\ref{S-4}, we get then the following:
\begin{Thm}
$(\ucalG, r,l,\sbullet,\omega_\calG)$ 
is a symplectic groupoid for $(U,\alpha)$.
\end{Thm}

\begin{Rem}
It is interesting to note that the above theorem holds also without
Assumption~\ref{ass-sim}, as can be proved directly.
However, $\calG$ as we have defined it is not the phase space for
$(U,\alpha)$ in the general case, the missing information being a class
of homotopic paths inside a symplectic leaf joining the given
endpoints.
\end{Rem}


\end{document}